\documentclass{amsart}
\usepackage{amsfonts}

\theoremstyle{plain}

\numberwithin{equation}{section}

\input{tcilatex}

\begin{document}
\author{V P Belavkin}
\address{Department of Mathematics, University of Nottingham, UK}
\title[Quantum Probabilities and Paradoxes]{Quantum Probabilities and
Paradoxes of the Quantum Century}
\date{Received May 5, 2000}
\email{vpb@maths.nott.ac.uk}
\urladdr{http://www.maths.nottingham.ac.uk/personal/vpb/inaugural/index.htm}
\keywords{Quantum probability, Quantum statistics, Quantum paradoxes,
Quantum measurement}
\dedicatory{Dedicated to Robin Hudson on his 60th birthday}
\thanks{ Published in: \textit{Inf. Dim. Anal., Quntum Probability and
Related Topics}, \textbf{3} (4) 577--610 (2000)}
\maketitle

\begin{abstract}
A history and drama of the development of quantum probability theory is
outlined starting from the discovery of the Plank's constant exactly a 100
years ago. It is shown that before the rise of quantum mechanics 75 years
ago, the quantum theory had appeared first in the form of the statistics of
quantum thermal noise and quantum spontaneous jumps which have never been
explained by quantum mechanics. Moreover, the only reasonable probabilistic
interpretation of quantum theory put forward by Max Born was in fact in
irreconcilable contradiction with traditional mechanical reality and
classical probabilistic causality. This led to numerous quantum paradoxes,
some of them due to the great inventors of quantum theory such as Einstein
and Schroedinger. They are reconsidered in this paper from the modern
quantum probabilistic point of view.
\end{abstract}

\tableofcontents

\section{Introduction: The Common Thread of Mathematical Sciences}

\medskip

\begin{quote}
\textit{The whole is more than the sum of its parts} -- Aristotle.
\end{quote}

This is the famous superadditivity law from Aristotle's \textit{Metaphysics}
which studies `the most general or abstract features of reality and the
principles that have universal validity'. Certainly in this broad definition
mathematics and physics are parts of metaphysics.

The aim of this lecture is to demonstrate, by the example of the history and
drama of the development of Quantum Theory during last century, that this
law is also applicable to metaphysics itself, as the unification of
different branches of mathematics and physics, which is more then just the
sum of pure and applied mathematics, statistics and mathematical physics.\
Quantum theory is a mathematical theory which studies the most fundamental
features of reality in a unified form of waves and matter it raises and
solves the most fundamental riddles of Nature by developing and utilizing
mathematical concepts and methods of all branches of modern mathematics,
including statistics.

Indeed, as we shall see, it began with the discovery of new laws for
``quantum'' numbers, the natural objects \ which are the foundation of pure
mathematics. (God made the integers; the rest is man's work). Next it
invented new applied mathematical methods for solving quantum mechanical
matrix and partial differential equations. Next it married probability with
algebra to obtain unified treatment of waves and particles in nature, giving
birth to quantum probability and creating new branches of mathematics such
as quantum logics, quantum topologies, quantum geometries, quantum groups.
It inspired the recent creation of quantum analysis and quantum calculus, as
well as quantum statistics and quantum stochastics.

Specialists in different narrow branches of mathematics rarely understand
quantum theory as a common thread which runs through everything. The
creators of quantum mechanics, the theory invented for interpretation of the
dynamical laws of fundamental particles, were unable to find a consistent
interpretation of it since they were physicists with a classical
mathematical education. After inventing quantum mechanics they spent much of
their lives trying to tackle the Problem of Quantum Measurement, the
greatest problem of quantum theory, not just of quantum mechanics, or even
of unified quantum field theory, which would be the same ``thing in itself''
as quantum mechanics of closed systems without such interpretation. As we
shall see, the solution to this problem can be found in the framework of
Quantum Probability as a part of a unified mathematics rather than physics.
Most modern mathematical physicists have a broad mathematical education, but
it ignores just two crucial aspects -- information theory and statistical
conditioning. So they gave up this problem as an unsolvable -- and it is
indeed unsolvable in the traditional framework of mathematical physics.

In order to appreciate the quantum drama which has been developing through
the whole century, and to estimate possible consequences of it in the new
quantum technological age, it seems useful to give a brief account of the
discovery of quantum theory at the beginning of the 20th century.

\section{The Quantum Century Begins}

\medskip

\begin{quote}
\textit{In science one tries to tell people, in such a way as to be
understood by everyone, something that no one even knew before. But in
poetry, it's quite opposite} -- Paul Dirac.
\end{quote}

Quantum Theory is the greatest intellectual achievement of the past century.
Since the discovery of quanta by Max Planck exactly 100 years ago \cite%
{Plnkqp} on the basis of spectral analysis of quantum thermal noise it has
produced numerous paradoxes and confusions even in the greatest scientific
minds such as Einstein, de Broglie, Schr\"{o}dinger, Bell, and it still
confuses many contemporary philosophers and scientists. The rapid
development of a beautiful and sophisticated mathematics for quantum
mechanics and the development of its interpretation by Bohr, Born,
Heisenberg, Dirac and many others who abandoned traditional causality, were
little help in resolving these paradoxes despite the astonishing success in
the prediction of quantum phenomena. Both the implication and consequences
of the quantum theory of light and matter, as well as its profound
mathematical, conceptual and philosophical foundations are not yet
understood completely. As Planck stated later:-

\begin{quote}
\textit{If anyone says he can think about quantum problems without getting
giddy, that only shows he has not understood the first thing about them}
\end{quote}

\subsection{ The Discovery of Matrix Mechanics}

In 1905 Einstein, examining the photoelectric effect, proposed a quantum
theory of light, only later realizing that Planck's theory made implicit use
of this quantum light hypothesis. Einstein saw that the energy changes in a
quantum material oscillator occur in jumps which are multiples of $\omega $.

In 1912 Niels Bohr worked in the Rutherford group in Manchester on his
theory of the electron in an atom. He was puzzled by the discrete spectra of
light which is emitted by atoms when they are subjected to an excitation. He
was influenced by the ideas of Planck and Einstein and addressed a certain
paradox in his work. How can energy be conserved when some energy changes
are continuous and some are discontinuous, i.e. change by quantum amounts?%
\textit{\ }Bohr conjectured that an atom could exist only in a discrete set
of stable energy states, the differences of which amount to the observed
energy quanta. Bohr returned to Copenhagen and published a revolutionary
paper on the hydrogen atom in the next year. He suggested his famous formula 
\begin{equation*}
E_{m}-E_{n}=\hbar \omega _{mn} 
\end{equation*}
from which he derived the major laws which describe physically observed
spectral lines. This work earned Niels Bohr the 1922 Nobel Prize about $%
10^{5}$ Swedish Kroner. This was a time when it was still possible to get
such a sum for so simple an equation!

In 1925 a young German theoretical physicist, Heisenberg, gave a preliminary
account of a new and highly original approach to the mechanics of the atom 
\cite{Heis25}. He was influenced by Niels Bohr and proposed to substitute
for the position coordinate of an electron an array 
\begin{equation*}
q_{mn}\left( t\right) =q_{mn}e^{i\omega _{mn}t}. 
\end{equation*}
His Professor, Max Born, was a mathematician who immediately recognized an
infinite matrix algebra in Heisenberg's arrays $Q=\left[ q_{mn}\right] $.
The classical momentum was also replaced by a similar matrix, 
\begin{equation*}
P\left( t\right) =\left[ p_{mn}e^{i\omega _{mn}t}\right] 
\end{equation*}
and $P$ and $Q$ matrices were postulated to follow a commutation law: 
\begin{equation*}
\left[ Q\left( t\right) ,P\left( t\right) \right] =i\hbar I, 
\end{equation*}
where $I$ is the unit matrix. The classical Hamiltonian equations of
dynamical evolution were now replaced by 
\begin{equation*}
\frac{d}{dt}Q\left( t\right) =\frac{i}{\hbar }\left[ H,Q\left( t\right) %
\right] ,\;\frac{d}{dt}P\left( t\right) =\frac{i}{\hbar }\left[ H,P\left(
t\right) \right] , 
\end{equation*}
where $H=\left[ E_{n}\delta _{mn}\right] $ is the diagonal Hamilton matrix 
\cite{BHJ26}. Thus quantum mechanics was first invented in the form of 
\textit{matrix mechanics}, emphasizing the possibilities of quantum
transitions, or jumps between the stable energy states $E_{n}$ of an
electron. In 1932 Heisenberg was awarded the Nobel Prize for his work in
mathematical physics.

Conceptually, the new atomic theory was based on the positivism of Mach as
it operated not with real space-time but with only observable quantities
like atomic transitions. However many leading physicists were greatly
troubled by the prospect of loosing reality and deterministic causality in
the emerging quantum physics. Einstein, in particular, worried about the
element of `chance' which had entered physics. In fact, this worries came
rather late since Rutherford had introduced a spontaneous effect when
discussing radio-active decay in 1900.

Thus, quantum theory first emerged as the result of experimental data not in
the form of quantum mechanics but in the form of statistical observations of
quantum noise, the basic concept of quantum probability and quantum
stochastic processes.

\subsection{The Discovery of Wave Mechanics}

The corpuscular nature of light seemed to contradict the Maxwell
electromagnetic wave theory of light. In 1924 Einstein wrote:

\begin{quote}
\textit{There are therefore now two theories of light, both indispensable,
and - as one must admit today, despite twenty years of tremendous effort on
the part of theoretical physics - without any logical connection.}
\end{quote}

The solution to the paradox of the wave/corpuscular duality of light came
unexpectedly when de Broglie made in 1923 the even more bizarre conjecture
of extending this duality also to material particles. He used the
Hamilton-Jacobi theory which had been applied both to particles and waves.

In 1925, Schr\"{o}dinger gave a seminar on de Broglie's material waves, and
a member of the audience suggested that there should be a wave equation.
Within a few weeks Schr\"{o}dinger found his celebrated wave equation, first
in a relativistic, and then in the non-relativistic form \cite{Schr26}.
Instead of seeking the classical solutions to the Hamilton-Jacobi equation 
\begin{equation*}
H\left( q,\frac{\hbar }{i}\frac{\partial }{\partial q}\ln \psi \right) =E 
\end{equation*}
he suggested finding those wave functions $\psi $ which satisfy the equation 
\begin{equation*}
H\left( q,\frac{\hbar }{i}\frac{\partial }{\partial q}\right) \psi =E\psi 
\end{equation*}
(It coincides with the former equation only if $H$ is linear with respect to 
$p$). He also obtained the non-stationary wave equation written in terms of
the Hamiltonian operator $H$ as 
\begin{equation*}
i\hbar \frac{\partial }{\partial t}\psi \left( t\right) =H\psi \left(
t\right) . 
\end{equation*}

Schr\"{o}dinger published his revolutionary \textit{wave mechanics} in a
series of six papers \cite{Schr26c} in 1926 during a short period of
sustained creative activity that is without parallel in the history of
science. Like Shakespeare, whose sonnets were inspired by a dark lady, Schr%
\"{o}dinger was inspired by a mysterious young lady of Arosa where he took
ski holidays during the Christmas 1925 but ``had been distracted by a few
calculations''. This was the second formulation of quantum theory, which he
successfully applied to the Hydrogen atom, oscillator and other quantum
mechanical systems. The mathematical equivalence between the two
formulations of quantum mechanics was understood by Schr\"{o}dinger
immediately, and he also introduced operators associated with each dynamical
variable.

As a young girl said later to Schr\"{o}dinger, who discovered the quantum
mechanics of wave matter:

\begin{quote}
\textit{Hey, you never even thought when you began that so much sensible
stuff would come out of it.}
\end{quote}

Unlike Heisenberg and Born's matrix mechanics, the general reaction towards
wave mechanics was immediately enthusiastic. Plank described Schr\"{o}%
dinger's wave mechanics as ``epoch-making work''. Einstein wrote:- ``the
idea of your work springs from true genius...''. Next year Schr\"{o}dinger
was nominated for the Nobel Prize, but he failed to receive it in this and
five further consecutive years of his nominations by most distinguished
physicists of the world, the reason behind his rejection being ``the highly
mathematical character of his work''. Only in 1933 did he receive his prize,
this time jointly with Dirac, and this was the first, and perhaps the last,
time when the Nobel Prize for physics was given to true mathematical
physicist.

\subsection{Interpretations of Quantum Mechanics}

The creators of the rival matrix quantum mechanics were forced to accept the
simplicity and beauty of Schr\"{o}dinger's approach. In 1926 Max Born put
forward the statistical interpretation of the wave function by introducing
the statistical mean 
\begin{equation*}
\left\langle H\right\rangle =\int \bar{\psi}\left( x\right) H\psi \left(
x\right) \mathrm{d}x 
\end{equation*}
for a dynamical variable $H$ in the physical state, described by $\psi $.
This was developed in Copenhagen and gradually was accepted by almost all
physicists as the ``Copenhagen interpretation''. Born by education was a
mathematician, and he would be the only mathematician ever to receive the
Nobel Prize, and then only in 1953, for his statistical studies of wave
functions, if he did not become a physicist, later Professor of Natural
Philosophy at Edinburgh. Bohr, Born and Heisenberg considered electrons and
quanta as unpredictable particles which cannot be visualized in the real
space and time.

Schr\"{o}dinger was a champion of the idea that the most fundamental laws of
the microscopic world are absolutely random before he discovered wave
mechanics. While he was preparing his Inaugural lecture for December 9, 1922
at the University of Z\"{u}rich, he wrote to Pauli:

\begin{quote}
\textit{I for my part believe, horribile dictu, that the energy -- momentum
law is violated in the process of radiation.}
\end{quote}

(quoted on page 152 in \cite{Schr})

Following de Broglie, Schr\"{o}dinger initially thought that the wave
function corresponds to a physical vibration process in a real continuous
space-time because it was not stochastic, but he was puzzled by the failure
to explain the blackbody radiation and photoelectric effect from this wave
point of view. In fact the wave interpretation applied to light quanta leads
back to classical electrodynamics, as his relativistic wave equation for a
single photon coincides mathematically with the classical wave equation.
However after realizing that the time-dependent $\psi $ is a complex
function in his fourth 1926 paper \cite{Schr26c}, submitted just a few days
before Born's, he admitted that the wave function $\psi $ cannot be given a
direct interpretation, and described the wave density $\psi \bar{\psi}$ as a
sort of weight function for superposition of point-mechanical configurations.

Bohr invited Schr\"{o}dinger to Copenhagen and tried to convince him of the
particle-probabilistic interpretation of quantum mechanics. The discussion
between them went on day and night, without reaching any agreement. The
conversation, however deeply affected both men. Schr\"{o}dinger recognized
the necessity of admitting both wave and particles, but he never devised a
comprehensive interpretation rival to Copenhagen orthodoxy. \ Bohr ventured
more deeply into philosophical waters and emerged with his concept of
complementarity:

\begin{quote}
\textit{Evidence obtained under different experimental conditions cannot be
comprehended within a single picture, but must be regarded as complementary
in the sense that only the totality of the phenomena exhausts the possible
information about the objects.}
\end{quote}

Later Schr\"{o}dinger accepted the probabilistic interpretation of $\psi 
\bar{\psi}$, although he did not consider these probabilities classically,
but as the strength of our belief or anticipation of an experimental result.
In this sense the probabilities are closer to propensities than to the
frequencies of the statistical interpretation of Born and Heisenberg. Schr%
\"{o}dinger had never accepted the subjective positivism of Bohr and
Heisenberg, and his philosophy is closer to that called representational
realism. He was content to remain a critical unbeliever. For a deeper
analysis of probabilistic roots in interpretation of quantum mechanics see 
\cite{Acc97}.

The most outspoken opponent of a/the probabilistic interpretation was
Einstein. Albert Einstein admired the new development of quantum theory but
with suspicion, and rejected its acausality and probabilistic
interpretation. It was against his scientific instinct to accept quantum
mechanics with its statistical interpretation as a complete description of
physical reality. There are famous sayings of him on that account:

\begin{quote}
\textquotedblleft \textit{...God is subtle but he is not
malicious\textquotedblright ...\textquotedblleft God doesn't play
dice\textquotedblright ...}
\end{quote}

In the famous debates on the probabilistic interpretation of quantum
mechanics of Einstein and Niels Bohr, Schr\"{o}dinger was often taking the
side of his friend Einstein, and this may explain why he was distancing
himself form the statistical interpretation of his wave function.

There have been many other attempts to retain the deterministic realism in
the quantum world, the most extravagant among these being the many world
interpretation. Certainly there are some advantages living in many worlds:
One can give, for example, many inaugural lectures if allowed to jump from
one world into another. From the philosophical and practical point of view,
however, to have an infinite (continuum?) number of real worlds at the same
time seems not better than to have none.

\section{Quantum Probabilities and Paradoxes}

\medskip

\begin{quote}
\textit{How wonderful we have met with a paradox, now we have some hope of
making progress} - Niels Bohr.
\end{quote}

In 1932 von Neumann put quantum theory on firm theoretical basis by setting
the mathematical foundation for new, quantum, probability theory, the
quantitative theory for counting non commuting events of quantum logics.
This noncommutative probability theory is based on essentially more general
axioms than the classical (Kolmogorovian) probability of commuting events,
which form common sense Boolean logic, first formalized by Aristotle. It has
been under extensive development during the last 30 years since the
introduction of algebraic and operational approaches for treatment of
noncommutative probabilities, and currently serves as the mathematical basis
for quantum information and measurement theory.

\subsection{Uncertainties and new logics}

Bohr was concerned with the paradox of spontaneous emission. He addressed
the question: How does the electron know when to emit radiation? Bohr, Born
and Heisenberg abandoned causality of traditional physics in the most
positivistic way. Max Born said:-

\begin{quote}
\textit{If God made the world a perfect mechanism, ... we need not solve
innumerable differential equations, but can use dice with fair success.}
\end{quote}

\subsubsection{Heisenberg uncertainty relations}

In 1927 Heisenberg derived \cite{Heis27} his famous uncertainty relation 
\begin{equation*}
\Delta Q\Delta P\geq \hbar /2,\quad \Delta T\Delta E\geq \hbar /2 
\end{equation*}
which gave mathematical support to the revolutionary complementary principle
of Bohr. The first relation was easily proved in the Schr\"{o}dinger
representations $Q=x$, $P=\frac{\hbar }{i}\frac{\partial }{\partial x}$ in
terms of the standard deviations 
\begin{equation*}
\Delta Q=\left( \left\langle Q^{2}\right\rangle -\left\langle Q\right\rangle
^{2}\right) ^{1/2},\quad \Delta P=\left( \left\langle P^{2}\right\rangle
-\left\langle P\right\rangle ^{2}\right) ^{1/2}. 
\end{equation*}
The second relation, which was first stated by analogy of $t$ with $x$ and
of $E$ with $P$, can be proved \cite{Be76, Hol80} in terms of the optimal
measurement of the initial time as an unknown parameter $\tau $ of the Schr%
\"{o}dinger's state $\psi \left( t-\tau \right) $ which is realized by the
measurement of self-adjoint operator $T=t$ in an extended representation
where the Hamiltonian $H$ is represented by the operator $E=i\hbar \frac{%
\partial }{\partial t}$. As Dirac stated:

\begin{quote}
\textit{Now when Heisenberg noticed that, he was really scared.}
\end{quote}

Einstein launched an attack on the uncertainty relation at the Solvay
Congress in 1927, and then again in 1930, by proposing cleverly devised
thought experiments which would violate this relation. Most of these
imaginary experiments were designed to show that interaction between the
microphysical object and the measuring instrument is not so inscrutable as
Heisenberg and Bohr maintained. He suggested, for example, a box filled with
radiation with a clock. The clock is designed to open a shutter and allow
one photon to escape. By weighing the box the photon energy and the time of
escape can both be measured with arbitrary accuracy.

After proposing this argument Einstein is reported to have spent a happy
evening, and Niels Bohr an unhappy one. After a sleepless night he showed
next morning that Einstein was not right. Mathematically his solution can be
expressed by the following formula 
\begin{equation*}
X=t+Q,\quad Y=i\hbar \frac{\partial }{\partial t}+P 
\end{equation*}
for the measuring quantity $X$, the pointer coordinate of the clock, and the
observable $Y$ for indirect measurement of photon energy $E=i\hbar \frac{%
\partial }{\partial t}$ in the Einstein experiment, where $Q$ and $P$ are
the position and momentum operators of the compensation weight under the
box. Due to the initial independence of the weight, the commuting
observables $X$ and $Y$ have even greater uncertainty 
\begin{equation*}
\Delta X\Delta Y\geq \Delta T\Delta E+\Delta Q\Delta P\geq \hbar 
\end{equation*}
than that predicted by Heisenberg uncertainty $\Delta T\Delta E\geq \hbar /2$%
.

\subsubsection{Nonexistence of hidden variables}

Einstein hoped that eventually it would be possible to explain the
uncertainty relations by expressing quantum mechanical observables as
functions of some hidden variables $\lambda $ in deterministic physical
states such that the statistical aspect will arise as in classical
statistical mechanics by averaging these observables over $\lambda $.

Von Neumann's monumental book \cite{Neum32} on the mathematical foundations
of quantum theory was therefore a timely contribution, clarifying, as it
did, this point. \ Inspired by Lev Landau, he introduced, for the unique
characterization of the statistics of a quantum ensemble, the statistical
density operator $\rho $ which eventually, under the name state, became a
major tool in quantum statistics. He considered the linear space $\mathcal{L}
$ of all bounded observables in a quantum system described by the Hermitian
operators $L^{\dagger }=L$ in a Hilbert space $\mathfrak{h}$, and defined
the expectation $\left\langle L\right\rangle $ of each $L\in \mathcal{L}$ in
a state $\rho $ by ultra-weakly continuous functional $\left\langle
L\right\rangle =\mathrm{Tr}L\rho $. Then he noted that any\emph{\ physically
continuous} additive functional $L\mapsto \left\langle L\right\rangle $ has
such trace form, and proved that in order to have positive probabilities $%
\Pr \left( E\right) $ for all quantum mechanical events as expectations $%
\left\langle E\right\rangle $ of yes-no observables $E\in \mathcal{L}$
described by Hermitian projectors $E^{2}=E$ and probability one for the
identity event described by identity operator $I$, i.e. 
\begin{equation*}
\Pr \left( E\right) \geq 0,\quad \Pr \left( I\right) =1, 
\end{equation*}
the statistical operator $\rho $ must be positive-definite and have trace
one. He applied this technique to the analysis of the completeness problem
of quantum theory, i.e. whether it constitutes a logically closed theory or
whether it could be reformulated as an entirely deterministic theory through
the introduction of hidden parameters (additional variables which, unlike
ordinary observables, are inaccessible to measurements). He came to the
conclusion that

\begin{quote}
\textit{the present system of quantum mechanics would have to be objectively
false, in order that another description of the elementary process than the
statistical one may be possible}

(quoted on page 325 in \cite{Neum32})
\end{quote}

To prove this theorem, von Neumann showed that there is no such state which
would be dispersion-free simultaneously for all quantum events $E\in 
\mathcal{L}$. For each such state, he argued, 
\begin{equation*}
\left\langle E^{2}\right\rangle =\left\langle E\right\rangle =\left\langle
E\right\rangle ^{2} 
\end{equation*}
for all events $E$ would imply that $\rho =0$, which cannot be statistical
operators as $0$ has trace $0$. Thus no state can be considered as a mixture
of dispersion-free states, each of them associated with a definite value of
hidden parameters. There are simply no such states, and thus, no hidden
parameters. In particular this implies that the statistical nature of pure
states, which are described by one-dimensional projectors $\rho =P_{\psi }$
corresponding to wave functions $\psi $, cannot be removed by supposing them
to be a mixture of dispersion-free substates.

It is widely believed that in 1966 John Bell showed that von Neuman's proof
was in error, or at least his analysis left the real question untouched \cite%
{Bell66}. To discredit the von Neumann's proof he constructed an example of
dispersion-free states parametrized for each quantum state $\rho $ by a real
parameter $\lambda $ for \ the case of two dimensional $\mathfrak{h}$. He
succeeded to do this by weakening the assumption of the additivity for such
states, requiring it only for the commuting observables in $\mathcal{L}$,
and by abandoning the linearity of the constructed expectations in $\rho $
described by spin polarization vector $\mathbf{r}$. There is no reason, he
argued, to keep the linearity in $\rho $ for the observable eigenvalues
determined by $\lambda $ and $\rho $, and to demand the additivity for
non-commuting observables as they are not simultaneously measurable: The
measured eigen-values of a sum of noncommuting observables are not the sums
of the eigen-values of this observables measured separately. Bell found for
each spin-projection $L=\sigma \left( \mathbf{e}\right) $ a family $%
s_{\lambda }\left( \mathbf{e}\right) $, $\left| \lambda \right| \leq 1/2$ of
dispersion free values $\pm 1$, reproducing the expectation $\left\langle
\sigma \left( \mathbf{e}\right) \right\rangle =\mathbf{e\cdot r}$ in the
pure quantum state when uniformly averaged over the $\lambda $. However his
example does not contradict the von Neumann theorem even if the latter is
strengthened by the restriction of the additivity only to the commuting
observables: The constructed dispersion-free expectation function $L\mapsto
\left\langle L\right\rangle _{\lambda }$ is not \emph{physically continuous}
on $\mathcal{L}$ because the value $\left\langle L\right\rangle _{\lambda
}=s_{\lambda }\left( \mathbf{e}\right) $ is one of the eigen-values $\pm 1$
for each $\lambda $, and it covers both values when the directional vector $%
\mathbf{e}$ rotates over the three-dimensional sphere. A function $\mathbf{%
e\mapsto }\left\langle \sigma \left( \mathbf{e}\right) \right\rangle
_{\lambda }$ on the continuous manifold (sphere) with discontinuous values
can be continuous only if it is constant, but this is ruled out by the
demand to reproduce the variable in $\mathbf{e}$ expectation $\left\langle
\sigma \left( \mathbf{e}\right) \right\rangle =\mathbf{e\cdot r}$ by
averaging the constant $\left\langle \sigma \left( \mathbf{e}\right)
\right\rangle _{\lambda }$ over all $\lambda $. Measurements of the
projections of spin on the physically close directions should be described
by close expected values in any physical state specified by $\lambda $,
otherwise it cannot have physical meaning!

Since then there were innumerable attempts to introduce hidden variables in
ever more sophisticated forms, perhaps not yet discovered, which would
determine the complementary variables if the hidden variables were measured
precisely. In higher dimensions of $\mathfrak{h}$ all these attempts are
ruled out by Gleason's theorem \cite{Glea57} who proved that there is not
even one additive zero-one probability function if $\dim \mathfrak{h}>2$.

\subsubsection{Complementarity and common sense}

In view of the decisive importance of this analysis for the foundations of
quantum theory, Birkhoff and von Neumann \cite{BiNe36} setup a system of
formal axioms for the calculus of logico-theoretical propositions concerning
results of possible measurements in a quantum system. They started with
formalizing the calculus of quantum events $E\in \mathcal{L}$ described by
orthoprojectors $E$, i.e. projective operators $E=E^{2}$ which are
orthogonal to the complements $E^{\bot }=I-E$ in the sense $E^{\dagger
}E^{\bot }=O$, that is $E^{\dagger }=E^{\dagger }E=E$ on the Hilbert space $%
\mathfrak{h}$. These are the Hermitian projectors $E\in \mathcal{E}$ which
have only two eigenvalues $\left\{ 0,1\right\} $ (``no'' and ``yes''). Such
calculus coincides with the calculus of the linear subspaces $\mathfrak{e}%
\subset \mathfrak{h}$ given by the ranges $\mathfrak{e}=E\mathfrak{h}$ of
the orthoprojectors $E$, in the same sense as the propositional calculus of
classical events coincides with the calculus of subsets of a Boolean
algebra. \ In this calculus the logical ordering $E\leq F$ implemented by
the algebraic relation $EF=E$ coincides with 
\begin{equation*}
\mathrm{range}\left( E\right) \subseteq \mathrm{range}\left( F\right) , 
\end{equation*}
the conjunction $E\wedge F$ corresponds to the intersection, 
\begin{equation*}
\mathrm{range}\left( E\wedge F\right) =\mathrm{range}\left( E\right) \cap 
\mathrm{range}\left( F\right) , 
\end{equation*}
however the disjunction $E\vee F$ is represented by the linear sum $%
\mathfrak{e}+\mathfrak{f}$ of the corresponding subspaces but not their
union, 
\begin{equation*}
\mathrm{range}\left( E\vee F\right) \neq \mathrm{range}\left( E\right) \cup 
\mathrm{range}\left( F\right) . 
\end{equation*}
Note that although $\mathfrak{e}+\mathfrak{f}=\mathrm{range}\left(
E+F\right) $, the operator $E+F$ does not coincide with orthoprojector $%
E\vee F$ onto $\mathfrak{e}+\mathfrak{f}$ unless $EF=0$. This implies that
the distributive law characteristic for propositional calculus of classical
logics no longer holds, but it still holds for compatible events, described
by commutative orthoprojectors due to the orthomodularity 
\begin{equation*}
E\leq I-F\leq G\Longrightarrow \left( E\vee F\right) \wedge G=E\vee \left(
F\wedge G\right) . 
\end{equation*}

Two events $E,F$ are called complementary if $E\vee F=I$, orthocomplementary
if $E+F=I$, incompatible or disjunctive if $E\wedge F=0$, and contradictory
or orthogonal if $EF=0$. As in the classical, common sense, logic
contradictory events are incompatible. However \emph{incompatible
propositions or events are not necessary contradictory }as can be easily
seen for any two nonorthogonal but not coinciding one-dimensional subspaces.
In particular, in quantum logics there exist complementary incompatible
pairs $E,F$, $E\vee F=I$, $E\wedge F=0$ which are not ortho-complementary in
the sense $E+F\neq I$, i.e. $EF\neq 0$ (it would be impossible in the
classical case). This is a rigorous logico-mathematical proof of Bohr's
complementarity.

As an example, we can consider the statement that a quantum system is in a
stable energy state $E$, and an incompatible proposition $F$, that it
collapses at a given time $t$, say. The incompatibility $E\wedge F=0$
follows from the fact that there is no state in which the system would
collapse preserving its energy, however these two propositions are not
contradictory (i.e. not orthogonal, $EF\neq 0$): the system might not
collapse if it is in other than $E$ stable state (remember the Schr\"{o}%
dinger's earlier belief that the energy law is valid only on average, and is
violated in the process of radiation).

In 1952 Wick, Wightman, and Wigner \cite{WWW52} showed that there are
physical systems for which not every orthoprojector corresponds to an
observable event, so that not every orthoprojector $P_{\psi }$ corresponding
to a wave function $\psi $ is a pure state. This is equivalent to the
admission of some selective events which are dispersion-free in all pure
states. Jauch and Piron \cite{JaPi63} incorporated this situation into
quantum logics and proved in the context of this most general approach that
the hidden variable interpretation is only possible if the theory is
observably wrong, i.e. if incompatible events are in fact compatible or
contradictory.

Bell criticized this also, as well as he criticized the Gleason's theorem,
but this time his arguments were not based even on the classical ground of
usual probability theory. Although he explicitly used the additivity of the
probability on the orthogonal events in his counterexample for $\mathfrak{h}=%
\mathbb{C}^{2}$, but he questioned : `That so much follows from such
apparently innocent assumptions leads us to question their innocence'. \cite%
{Bell66}. In fact this was equivalent to questioning of the additivity of
the classical probability on the corresponding disjoint subsets, but he
didn't suggest any other complete system of physically reasonable axioms for
introducing such peculiar ``nonclassical'' hidden variables, not even a
single counterexample to the orthogonal nonadditivity even for the simplest
case $\dim \mathfrak{h}=2$. Thus Bell implicitly rejected classical
probability theory in the quantum world, but he didn't want to accept
quantum probability as the only possible theory for explaining the
microworld. Obviously that even if such attempt was successful for a single
quantum system (as he thought his unphysical discontinuous construction in
the case $\mathfrak{h}=\mathbb{C}^{2}$ was), the classical composition low
of two such systems would not allow to extend the dispersion-free
product-states to the whole composed quantum system because of their
nonadditivity on the whole observables space $\mathcal{L}$. The quantum
composition law together with the orthoadditivity excludes the hidden
variable possibility even for $\mathfrak{h}=\mathbb{C}^{2}$, otherwise it
would contradict to Gleason's theorem for the case $\dim \left( \mathfrak{h}%
\otimes \mathfrak{h}\right) =4$), not even to mention a hidden variable
reproduction of nonseparable, entangled states). In fact it justifies the
von Neumann's additivity assumption for the states on the whole operator
algebra $\mathcal{B}\left( \mathfrak{h}\right) $.

\subsection{Quantum measurement and decoherence}

Heisenberg derived from the uncertainty relation that `the nonvalidity of
rigorous causality is necessary and not just consistently possible'. Max
Born even stated:-

\begin{quote}
\textit{One does not get an answer to the question, what is the state after
collision? but only to the question, how probable is a given effect of the
collision?}
\end{quote}

\subsubsection{Spooky action at distance}

After his defeat on uncertainty relations Einstein seemed to have become
resigned to statistical interpretation of quantum theory, and at the 1933
Solvay Congress he listened to Bohr's paper on completeness of quantum
theory without objections. Then, in 1935, he launched a brilliant and subtle
new attack in a paper \cite{EPR} with two young co-authors, Podolski and
Rosen, which is known as the EPR paradox that has become of major importance
to the world view of physics. They stated the following requirement for a
complete theory as a seemingly necessary one:

\begin{quote}
\textit{Every element of physical reality must have a counterpart in the
physical theory.}
\end{quote}

The question of completeness is thus easily answered as soon as soon as we
are able to decide what are the elements of the physical reality. \ EPR then
proposed a sufficient condition for an element of physical reality:

\begin{quote}
\textit{If, without in any way disturbing the system, we can predict with
certainty the value of a physical quantity, then there exists an element of
physical reality corresponding to this quantity.}
\end{quote}

Then they designed a thought experiment the essence of which is that two
quantum ``bits'', particles spins of two electrons say, are brought together
to interact, and after separation an experiment is made to measure the spin
orientation of one of them. The state after interaction is such that the
measurement result $\tau =\pm \frac{1}{2}$ of one particle uniquely
determination the spin $z$-orientation $\sigma =\mp \frac{1}{2}$ of the
other particle. EPR apply their criterion of local reality: since the value
of $\sigma $ can be predicted by measuring $\tau $ without in any way
disturbing $\sigma $, it must correspond to an existing element of physical
reality. Yet the conclusion contradicts a fundamental postulate of quantum
mechanics, according to which the sign of spin is not an intrinsic property
of a complete description of the spin but is evoked only by a process of
measurement. Therefore, EPR conclude, quantum mechanics must be incomplete,
there must be hidden variables not yet discovered, which determine the spin
as an intrinsic property. It seams Einstein was unaware of the von Neumann's
theorem, although they both had the positions at the same Institute for
Advanced Studies at Princeton being among the original six mathematics
professors appointed there in 1933.

Bohr carefully replied to this challenge by rejecting the assumption of
local physical realism as stated by EPR \cite{Bohr35}: `There is no question
of a mechanical disturbance of the system under investigation during the
last critical stage of the measuring procedure. But even at this stage there
is essentially a question of \textit{an influence on the very conditions
which define the possible types of predictions regarding the future behavior
of the system}'. This influence became notoriously famous as Bohr's \textit{%
spooky action at a distance}. He had obviously meant the semi-classical
model of measurement, when one can statistically infer the state of one
(quantum) part of a system immediately after observing the other (classical)
part, whatever the distance between them. In fact, there is no paradox of
``spooky action at distance'' in the classical case, the statistical
inference, playing the role of such immediate action, is simply based on the
Bayesian selection rule of a posterior state from the prior mixture of all
such states, corresponding to the possible results of the measurement. Bohr
always emphasized that one must treat the measuring instrument classically
(the measured spin, or another bit interacting with this spin, as a
classical bit), although the classical-quantum interaction should be
regarded as purely quantum. The latter follows from non-existence of
semi-classical Poisson bracket (i.e. classical-quantum potential
interaction). Schr\"{o}dinger clarified this point more precisely then Bohr,
and he followed in fact the mathematical pattern of von Neumann measurement
theory. EPR paradox is related to so called Bell inequality the
probabilistic roots of which was evidentiated in \cite{Acc81}.

\subsubsection{Releasing Schr\"{o}dinger's cat}

Motivated by EPR paper, in 1935 Schr\"{o}dinger published a three part essay 
\cite{Schr35} on `The Present Situation in Quantum mechanics'. He turns to
EPR paradox and analyses completeness of the description by the wave
function for the entangled parts of the system. (The word \emph{entangled}
was introduced by Schr\"{o}dinger for the description of nonseparable
states.) He notes that if one has pure states $\psi \left( \sigma \right) $
and $\varphi \left( \tau \right) $ for each of two completely separated
bodies, one has maximal knowledge, $\chi \left( \sigma ,\tau \right) =\psi
\left( \sigma \right) \varphi \left( \tau \right) $, for two taken together.
But the converse is not true for the entangled bodies, described by a
non-separable wave function $\chi \left( \sigma ,\tau \right) \neq \psi
\left( \sigma \right) \varphi \left( \tau \right) $:

\begin{quote}
\textit{Maximal knowledge of a total system does not necessary imply maximal
knowledge of all its parts, not even when these are completely separated one
from another, and at the time can not influence one another at all.}
\end{quote}

To make absurdity of the EPR argument even more evident he constructed his
famous burlesque example in quite a sardonic style. A cat is shut up in a
steel chamber equipped with a camera, with an atomic mechanism in a pure
state $\rho _{0}=P_{\psi }$ which triggers the release of a phial of cyanide
if an atom disintegrates spontaneously (it is assumed that it might not
disintegrate in a course of an hour with probability $\mathrm{Tr}\left(
EP_{\psi }\right) =1/2$). If the cyanide is released, the cat dies, if not,
the cat lives. Because the entire system is regarded as quantum and closed,
after one hour, without looking into the camera, one can say that the entire
system is still in a pure state in which the living and the dead cat are
smeared out in equal parts.

Schr\"{o}dinger resolves this paradox by noting that the cat is a
macroscopic object, the states of which (alive or dead) could be told apart
by a macroscopic observation are distinct from each other whether observed
or not. He calls this `the principle of state distinction' for macroscopic
objects, which is in fact the postulate that the directly measurable system
(consisting of cat) must be classical:

\begin{quote}
\textit{It is typical in such a case that an uncertainty initially
restricted to an atomic domain has become transformed into a macroscopic
uncertainty which can be resolved through direct observation.}
\end{quote}

The dynamical problem of the transformation of the atomic, or ``coherent''
uncertainty, corresponding to a probability amplitude $\psi \left( \sigma
\right) $, into a macroscopic uncertainty, corresponding to a mixed state $%
\rho $, is called quantum \emph{decoherence }problem$.$ In order to make
this idea clear, let us give the solution of the Schr\"{o}dinger's
elementary decoherence problem in the purely mathematical way. Instead of
the values $\pm 1/2$ for the spin-variables $\sigma $ and $\tau $ we shall
use the values $\left\{ 0,1\right\} $ corresponding to the states of a
classical ``bit'', the simplest nontrivial system in classical probability
or information theory.

Consider the atomic mechanism as a quantum ``bit'', $\mathfrak{h}=\mathbb{C}%
^{2}$, the pure states of which are described by $\psi $-functions of the
variable $\sigma \in \left\{ 0,1\right\} $ (if atom is disintegrated, $%
\sigma =1$ , if not, $\sigma =0$) with scalar (complex) values $\psi \left(
\sigma \right) $ defining the probabilities $\left| \psi \left( \sigma
\right) \right| ^{2}$ of the quantum elementary propositions corresponding
to $\sigma =0,1$. The Schr\"{o}dinger's cat is a classical bit with only two
pure states $\tau \in \left\{ 0,1\right\} $ which can be identified with the
probability distributions $\delta _{0}\left( \tau \right) $ when alive $%
\left( \tau =0\right) $ and $\delta _{1}\left( \tau \right) $ when dead $%
\left( \tau =1\right) $. These and other (mixed) states can also be
described by the complex amplitudes $\varphi \left( \tau \right) $, however
they are uniquely defined by the probabilities $\left| \varphi \left( \tau
\right) \right| ^{2} $ up to a phase function of $\tau $, the phase
multiplier of $\varphi $ commuting with all cat observables $c\left( \tau
\right) $, not just up to a phase constant as in the case of the atom (only
constants commute with all atomic observables $A\in \mathcal{L}\left( 
\mathfrak{h}\right) $). Initially the cat is alive, so its amplitude
(uniquely defined up to the phase factor by the square root of probability
distribution $\delta _{0}$ on $\left\{ 0,1\right\} $) is $\delta \left( \tau
\right) $ that is $1$ if $\tau =0$, and $0$ if $\tau =1$.

The only meaningful classical-quantum reversible interaction affecting not
atom but cat as it is said after the hour, is described as both bits were
quantum, by 
\begin{equation*}
U\left[ \psi \otimes \varphi \right] \left( \sigma ,\tau \right) =\psi
\left( \sigma \right) \varphi \left( \left( \tau \bigtriangleup \sigma
\right) \right) , 
\end{equation*}
where $\tau \bigtriangleup \sigma =\left| \tau -\sigma \right| =\sigma
\bigtriangleup \tau $ is the difference ($\func{mod}2$) on $\left\{
0,1\right\} $. Applied to the initial product-state $\psi \otimes \delta $
it has the resulting probability amplitude 
\begin{equation*}
\chi \left( \sigma ,\tau \right) =\psi \left( \sigma \right) \delta \left(
\tau \bigtriangleup \sigma \right) =0\quad \mathrm{if}\quad \sigma \neq \tau
. 
\end{equation*}
Despite the fact that the initial state was pure, $\chi _{0}=\psi \otimes
\delta $ corresponding to the Cartesian product $\left( \psi ,0\right) $ of
the initial pure states $\psi \in \mathfrak{h}$ and $\tau =0$, the
reversible unitary evolution $U$ induces the mixed state for the
quantum-classical system ``atom+cat'' described by the wave function $\chi
\in \mathfrak{h}\otimes \mathbb{C}^{2}$.

Indeed, the observables of such a system are operator-functions $X$ of $\tau 
$ with values $X\left( \tau \right) $ in $\sigma $-matrices, represented as
block-diagonal $\left( \sigma ,\tau \right) $-matrices $\hat{X}=\left[
X\left( \tau \right) \delta _{\tau ^{\prime }}^{\tau }\right] $ of the
multiplication $X\left( \tau \right) \chi \left( \cdot ,\tau \right) $ at
each point $\tau \in \left\{ 0,1\right\} $. This means that the amplitude $%
\chi $ induces the same expectations 
\begin{equation*}
\left\langle \hat{X}\right\rangle =\sum_{\tau }\chi \left( \tau \right)
^{\dagger }X\left( \tau \right) \chi \left( \tau \right) =\sum_{\tau }%
\mathrm{Tr}X\left( \tau \right) \varrho \left( \tau \right) =\mathrm{Tr}\hat{%
X}\hat{\varrho} 
\end{equation*}
as the block-diagonal density matrix $\hat{\varrho}=\left[ \varrho \left(
\tau \right) \delta _{\tau ^{\prime }}^{\tau }\right] $ of the
multiplication by 
\begin{equation*}
\varrho \left( \tau \right) =E\left( \tau \right) P_{\psi }E\left( \tau
\right) =\Pr \left( \tau \right) P_{E\left( \tau \right) \psi } 
\end{equation*}
where $\Pr \left( \tau \right) =\left| \psi \left( \tau \right) \right| ^{2}$%
, $E\left( \tau \right) =P_{\delta _{\tau }}$ is the projection operator 
\begin{equation*}
\left[ E\left( \tau \right) \psi \right] \left( \sigma \right) =\delta
\left( \tau \bigtriangleup \sigma \right) \psi \left( \sigma \right) =\psi
\left( \tau \right) \delta _{\tau }\left( \sigma \right) , 
\end{equation*}
and $P_{E\left( \tau \right) \psi }=P_{\delta _{t}}$ is also the projector
onto $\delta _{\tau }\left( \cdot \right) =\delta \left( \cdot
\bigtriangleup \tau \right) $. The $4\times 4$-matrix $\hat{\varrho}$ is a
mixture of two orthogonal projectors $P_{\tau }\otimes P_{\tau }$, $\tau
=0,1 $: 
\begin{equation*}
\hat{\varrho}=\left[ P_{\delta _{\tau }}\delta _{\tau ^{\prime }}^{\tau }\Pr
\left( \tau \right) \right] =\sum_{\tau =0}^{1}\Pr \left( \tau \right)
P_{\delta _{\tau }}\otimes P_{\delta _{\tau }}. 
\end{equation*}

\subsubsection{Von Neumann's projection postulate}

Inspired by Bohr's complementarity principle, von Neumann proposed even
earlier the idea that every quantum measuring process involves an
unanalysable element. He postulated \cite{Neum32} that, in addition to the
continuous causal propagation of the wave function generated by the Schr\"{o}%
dinger equation, during a measurement, due to an action of the observer on
the object, the function undergoes a discontinuous, irreversible
instantaneous change. Thus, just prior to a measurement of the event $F$,
disintegration of atom, say, the quantum pure state $P_{\psi }$ changes to
the mixed one 
\begin{equation*}
\rho =\lambda P_{E\psi }+\mu P_{F\psi }=E\rho E+F\rho F, 
\end{equation*}
where $E=I-F$ is the orthocomplement event, and $\lambda =\mathrm{Tr}E\rho $%
, $\mu =\mathrm{Tr}F\rho $ are the probabilities of $E$ and $F$. Such change
is projective as it shows the second part of this equation, and it is called
the von Neumann projection postulate.

This linear irreversible decoherence process should be completed by the
nonlinear, acausal random jump to one of the pure states 
\begin{equation*}
\rho \mapsto P_{E\psi }\text{, or }\rho \mapsto P_{F\psi } 
\end{equation*}
depending on whether the tested event $F$ is false (the cat is alive, $\psi
_{0}=\lambda ^{-1/2}E\psi $), or true (the cat is dead, $\psi _{1}=\mu
^{-1/2}F\psi $). Such last step is the posterior prediction, called \emph{%
filtering} of the decoherent mixture of $\psi _{0}$ and $\psi _{1}$ by
selection of only one result of the measurement, is an unavoidable element
in every measurement process relating the state of the pointer of the
measurement (in this case the cat) to the state of the whole system. This
assures that the same result would be obtained in case of immediate
subsequent measurement of the same event $F$. The resulting change of the
initial wave-function $\psi $ is described up to normalization by one of the
projections 
\begin{equation*}
\psi \mapsto E\psi ,\quad \psi \mapsto F\psi 
\end{equation*}
and is sometimes called the L\"{u}ders projection postulate \cite{Lud51}.

Although unobjectionable from the purely logical point of view the von
Neumann theory of measurement soon became the target of severe criticisms.
Firstly it seams radically subjective, postulating the spooky action at
distance (the filtering) in a purely quantum system instead of deriving it.
Secondly the analysis is applicable to only the idealized situation of
discrete instantaneous measurements.

However as we already mentioned when discussing the EPR paradox, the process
of filtering is free from conceptual difficulty if it is understood as the
statistical inference about a mixed state in an extended stochastic
representation of the quantum system as a part of a semiclassical one, based
upon the results of observation in its classical part. In order to
demonstrate this, we can return to the dynamical model of Schr\"{o}dinger's
cat, identifying the quantum system in question with the Schr\"{o}dinger's
atom. The event $E$ (the atom exists) corresponds then to $\tau =0$ (the cat
is alive), $E=E\left( 0\right) $, and the complementary event is $F=E\left(
1\right) $. This model explains that the origin of the von Neumann
irreversible decoherence $P_{\psi }\mapsto \rho $ of the atomic state is in
the ignorance of the result of the measurement described by the partial
tracing over the cat's Hilbert space $\mathbb{C}^{2}$: 
\begin{equation*}
\rho =\mathrm{Tr}_{\mathbb{C}^{2}}\hat{\varrho}=\sum_{\tau =0}^{1}\Pr \left(
\tau \right) P_{\delta _{\tau }}=\varrho \left( 0\right) +\varrho \left(
1\right) , 
\end{equation*}
where $\varrho \left( \tau \right) =\left| \psi \left( \tau \right) \right|
^{2}P_{\delta _{\tau }}$. It has the entropy $S\left( \rho \right) =\mathrm{%
Tr}\rho \log \rho ^{-1}$ of the compound state $\hat{\varrho}$ of the
combined semi-classical system prepared for the indirect measurement of the
disintegration of atom by means of cat's death: 
\begin{equation*}
S\left( \rho \right) =-\sum_{\tau =0}^{1}\left| \psi \left( \tau \right)
\right| ^{2}\log \left| \psi \left( \tau \right) \right| ^{2}=S\left( \hat{%
\varrho}\right) 
\end{equation*}
It is the initial coherent uncertainty in the pure quantum state of the atom
described by the wave-function $\psi $ which is equal to one bit if
initially $\left| \psi \left( 0\right) \right| ^{2}=1/2=\left| \psi \left(
1\right) \right| ^{2}$.

This dynamical model of the measurement which is due to von Neumann, also
interprets the filtering $\rho \mapsto \rho _{\tau }$ simply as the
conditioning 
\begin{equation*}
\rho _{\tau }=\varrho \left( \tau \right) /\Pr \left( \tau \right)
=P_{\delta _{\tau }} 
\end{equation*}
of the joint classical-quantum state $\varrho \left( \cdot \right) $ by the
Bayes formula which is applicable due to the commutativity of actually
measured observable (the life of cat) with any other observable of the
combined semi-classical system.

Thus the atomic decoherence is derived from the unitary interaction of the
quantum atom with the classical cat. The spooky action at distance,
affecting the atomic state by measuring $\tau $, is simply the result of the
statistical inference (prediction after the measurement) of the atomic
posterior state $\rho _{\tau }=P_{\delta _{\tau }}$: the atom disintegrated
if and only if the cat is dead. A formal derivation of the von-Neumann-L\"{u}%
ders projection postulate and the decoherence in the general case by
explicit construction of unitary transformation in the extended
semi-classical system in given in \cite{Be94, StBe96}.

\section{Causality and Prediction of Future}

\medskip \medskip

\textit{In mathematics you don't understand things. You just get used to them%
} - John von Neumann.

The quantum probability approach resolves the famous paradoxes of quantum
measurement theory in a constructive way by giving exact nontrivial models
for the statistical analysis of the quantum observation processes underlying
these paradoxes. Conceptually it is based upon a new idea of quantum
causality called the Nondemolition Principle which divides the world into
the classical past, forming the consistent histories, and the quantum
future, the state of which is predictable for each such history. The
differential analysis of these models is based on It\^{o} stochastic
calculus. The most recent mathematical development of these methods leads to
a profound quantum filtering and control theory in quantum open systems
which has found numerous applications in quantum statistics, optics and
spectroscopy, and is an appropriate tool for the solution of the dynamical
decoherence problem for quantum communications and computations.

\subsection{Reality and nondemolition principle}

Schr\"{o}dinger like Einstein was deeply concerned with the loss of reality
and causality in the positivistic treatment of quantum measuring process by
Heisenberg and Born. Schr\"{o}dinger's remained unhappy with Bohr's reply to
the EPR paradox, Schr\"{o}dinger's own analysis was:

\begin{quote}
\textit{It is pretty clear, if reality does not determine the measured
value, at least the measurable value determines reality.}
\end{quote}

\subsubsection{Compatibility and time arrow}

Von Neumann's projection postulate and its dynamical realization can be
generalized to include cases with continuous spectrum of values. In fact
there many such developments, we will only mention here the most general
operational approach to quantum measurements of Ludwig \cite{Lud68}, and its
mathematical implementation by Davies and Lewis \cite{DaLe70} in the
``instrumental'' form. The stochastic realization of the corresponding
completely positive reduction map $\rho \mapsto \varrho \left( \cdot \right) 
$, resolving the corresponding instantaneous quantum measurement problem,
can be found in \cite{Oz84, Be94}. Because of the crucial importance of
these realizations for developing understanding of the mathematical
structure and interpretation of modern quantum theory, analyze the
mathematical consequences which can be drawn from such schemes we need to.

The generalized reduction of the wave function $\psi \left( x\right) $,
corresponding to a complete measurement with discrete or continuous data $y$%
, is described by a function $V\left( y\right) $ whose values are linear
operators $\psi \mapsto V\left( y\right) \psi $ but not isometric, $V\left(
y\right) ^{\dagger }V\left( y\right) \neq I$, with the following
normalization condition. The resulting wave-function 
\begin{equation*}
\chi \left( x,y\right) =\left[ V\left( y\right) \psi \right] \left( x\right) 
\end{equation*}
as a is normalized with respect to a given measure $\mu $ on $y$ in the
sense 
\begin{equation*}
\iint \left| \chi \left( x,y\right) \right| ^{2}\mathrm{d}\mu \mathrm{d}%
\lambda =\int \left| \psi \left( x\right) \right| ^{2}\mathrm{d}\lambda 
\end{equation*}
for any probability amplitude $\psi $ (normalized with respect to a measure $%
\lambda $). In the discrete case such as the case of two-point variables $%
y=\tau $ (EPR paradox, or Schr\"{o}dinger cat with the projection-valued $%
V\left( \tau \right) =E\left( \tau \right) $) this is usually satisfied by
taking $\mu $ to be the counting measure: 
\begin{equation*}
\int_{y}V\left( y\right) ^{\dagger }V\left( y\right) \mathrm{d}\mu =I,\quad 
\mathrm{or}\quad \sum_{y}V\left( y\right) ^{\dagger }V\left( y\right) =I. 
\end{equation*}
As in that simple case the realization can be always constructed \cite{Be94}
in terms of a unitary transformation $U$ and a normalized wave function $%
\varphi _{0}$ such that 
\begin{equation*}
U\left[ \psi \otimes \varphi _{0}\right] \left( x,y\right) =\chi \left(
x,y\right) 
\end{equation*}
for any $\psi $. The additional system described by ``the pointer coordinate 
$y$ of the measurement apparatus'' can be regarded as classical (like the
cat) as the actual observables in question are the measurable functions $%
g\left( y\right) $ represented by commuting operators $\hat{g}$ of
multiplication by these functions. Note that such observables, extended to
the quantum part as $I\otimes \hat{g}$, are compatible with any possible
(future) event, represented by an orthoprojector $F\otimes I$. The
probabilities (or, it is better to say, the propensities) of all such events
are the same in all states whether an observable $\hat{g}$ was measured but
the result is not read, or it was not measured at all. In this sense the
measurements of $\hat{g}$ are called \textit{nondemolition} with respect to
the future observables $F$, they do not demolish the picture of the
possibilities, or propensities of $F$. But they are not necessary compatible
with of the initial operators $F\otimes I$ of the quantum system under the
question in the present representation $U\left( F\otimes I\right) U^{\ast }$
corresponding to the actual states $\chi =U\left( \psi \otimes \varphi
_{0}\right) $.

Indeed, the Heisenberg operators 
\begin{equation*}
G=U^{\ast }\left( I\otimes \hat{g}\right) U 
\end{equation*}
of the nondemolition observables in general do not commute with the past
operators $F\otimes I$ on the initial states $\chi _{0}=\psi \otimes \varphi
_{0}$. One can see this from the example of the Schr\"{o}dinger cat. The
``cat observables'' in Heisenberg picture are represented by commuting
operators $G=\left[ g\left( \sigma +\tau \right) \delta _{\sigma ^{\prime
}}^{\sigma }\delta _{\tau ^{\prime }}^{\tau }\right] $ of multiplication by $%
g\left( \sigma +\tau \right) $, where the sum $\sigma +\tau =\left| \sigma
-\tau \right| $ is modulo 2. They do not commute with $F\otimes I$ unless $F$
is also a diagonal operator $\hat{f}$ of multiplication by a function $%
f\left( \sigma \right) $ in which case 
\begin{equation*}
\left[ F,G\right] \chi _{0}\left( \sigma ,\tau \right) =\left[ f\left(
\sigma \right) ,g\left( \sigma +\tau \right) \right] \chi _{0}\left( \sigma
,\tau \right) =0,\quad \forall \chi _{0}\text{.} 
\end{equation*}
However the restriction of the possibilities in a quantum system to only the
diagonal operators $F=\hat{f}$ which would eliminate the time arrow in the
nondemolition condition amounts to the redundancy of the quantum
consideration as all such (possible and actual) observables can be
simultaneously represented as the functions of $\left( \sigma ,\tau \right) $
as in the classical case.

\subsubsection{Transition from possible to actual}

The analysis above shows that as soon as dynamics is taking into
consideration even in the form of just a single unitary transformation, the
measurement process needs to specify the arrow of time, what is the
predictable future and what is the reduced past, what is possible and what
is actual with respect to this measurement. As soon as a measured observable 
$Y$ is specified, all other operators which do not commute with $Y$ become
redundant as possible in future observables. The algebra of such potential
observables (not the state which stays invariant in the Schr\"{o}dinger
picture unless the selection due to an inference has taken place!) reduces
to the subalgebra commuting with $Y$, and \emph{this reduction doesn't
change the reality} (the wave function remains the same and induces the same
but now mixed state on the smaller - reduced algebra!). Possible observables
in an individual system are only those which are compatible with the actual
observables. This is another formulation of Bohr's complementarity. More
specifically this can be rephrased in the form of a dynamical postulate of
quantum causality called the Nondemolition Principle \cite{Be94} which we
first formulate for a single instant of time $t$ in a quite obvious form:

\begin{quote}
In the interaction representation of a quantum system by an algebra $%
\mathcal{A}$ of (necessarily not all) operators on a Hilbert space $%
\mathfrak{h}$ of the system with a measurement apparatus, causal, or
nondemolition observables are represented only by those operators $Y$ on $%
\mathfrak{h}$ which are compatible with $\mathcal{A}$: 
\begin{equation*}
\left[ X,Y\right] :=XY-YX=0,\quad \forall X\in \mathcal{A} 
\end{equation*}
(this is usually written as $Y\in \mathcal{A}^{\prime }$, where $\mathcal{A}%
^{\prime }$, called the commutant of $\mathcal{A}$, in this formulation is
not necessarily contained in $\mathcal{A}$). Each measurement process of the
history for a quantum system $\mathcal{A}$ can be represented as
nondemolition by the causal observables in the appropriate representation of 
$\mathcal{A}$ on a Hilbert space $\mathfrak{h}$.
\end{quote}

Note that the space of interaction representation $\mathfrak{h}$ plays here
the crucial role: the reduced operators 
\begin{equation*}
X_{0}=\left( I\otimes \varphi _{0}\right) ^{\ast }X\left( I\otimes \varphi
_{0}\right) ,\;Y_{0}=\left( I\otimes \varphi _{0}\right) ^{\ast }Y\left(
I\otimes \varphi _{0}\right) 
\end{equation*}
for commuting $X$ and $Y$ might not commute on the smaller space $\mathfrak{h%
}_{0}$ of the initial states $\psi \otimes \varphi _{0}$ with a fixed $%
\varphi _{0}\in \mathfrak{f}$. Even if the nondemolition observables $Y$ is
faithfully represented by $Y_{0}$ on initial space $\mathfrak{h}_{0}$, as it
is in the case $Y=G$ of the Schr\"{o}dinger's cat with $\varphi _{0}\left(
\tau \right) =\delta \left( \tau \right) $ where $Y_{0}$ is the
multiplication operator $G_{0}=\hat{g}$ for $\psi $: 
\begin{equation*}
G\left( \psi \otimes \delta \right) \left( \sigma ,\tau \right) =g\left(
\sigma +\tau \right) \psi \left( \sigma \right) \delta \left( \tau \right)
=g\left( \sigma \right) \psi \left( \sigma \right) \delta \left( \tau
\right) =\left( G_{0}\psi \otimes \delta \right) \left( \sigma ,\tau \right)
, 
\end{equation*}
there is no usually room in $\mathfrak{h}_{0}$ to represent all Heisenberg
operators $X\in \mathcal{A}$ commuting with $Y$ on $\mathfrak{h}$. The
induced operators $Y_{0}$ do not commute with all operators $F$ of the
system initially represented on $\mathfrak{h}_{0}$, and this is why the
measurement of $Y_{0}$ is thought to be demolition on $\mathfrak{h}_{0}$.
However in all such cases the future operators $X$ reduced to $X_{0}$ on $%
\mathfrak{h}_{0}$,\ commute with $Y_{0}$ as they are decomposable with
respect to $Y_{0}$, although the reduction $X\mapsto X_{0}$ is not the
Heisenberg but irreversible dynamical map. It can be explicitly seen for the
atom described by the Heisenberg operators $X=U^{\ast }\left( F\otimes
I\right) U$ in the interaction representation with the cat: 
\begin{equation*}
X_{0}=\sum_{\tau }E\left( \tau \right) FE\left( \tau \right) ,\quad
Y_{0}=\sum g\left( \tau \right) E\left( \tau \right) . 
\end{equation*}
Thus the nondemolition principle should not be considered as a restriction
on the possible observations for a given dynamics but rather a condition for
the causal dynamics to be compatible with the given observations. As was
proved in \cite{Be92a}, the causality condition is necessary and sufficient
for the existence of a conditional expectation for any state on the total
algebra $\mathcal{A}\vee \mathcal{B}$ with respect to a commutative
subalgebra $\mathcal{B}$ of nondemolition observables $Y$. Thus the
nondemolition causality condition amounts exactly to the existence of the
conditional states, i.e. to the predictability of the states on the algebra $%
\mathcal{A}$ upon the measurement results of the observables in $\mathcal{B} 
$. Then the transition from a prior $\rho $ to a posterior state $\rho
_{y}=P_{V\left( y\right) \psi }$ is simply the result of gaining a knowledge 
$y$ defining the actual state in the decoherent mixture 
\begin{equation*}
\rho =\int V\left( y\right) P_{\psi }V\left( y\right) ^{\ast }\mathrm{d}\mu
=\int P_{V\left( y\right) \psi }f\left( y\right) \mathrm{d}\mu , 
\end{equation*}
of all possible states, where $f\left( y\right) =\left\| V\left( y\right)
\psi \right\| ^{2}$ is the probability density of all possible $y$ defining
the output measure $\mathrm{d}\nu =f\mathrm{d}\mu $. As Heisenberg always
emphasized, ``quantum jump'' is contained in the transitions from possible
to actual.

If an algebra $\mathcal{B}$ of actual observables is specified at a time $t$%
, there must be a causal representation $\mathcal{B}_{t}$ of $\mathcal{B}$
on $\mathfrak{h}$ with respect to the present $\mathcal{A}_{t}$ and all
future possible representations $\mathcal{A}_{s}$, $s>t$ of the quantum
system on the same Hilbert space $\mathfrak{h}$ (they might not coincide
with $\mathcal{A}_{t}$ if the system is open \cite{AFL82}). The past
representations $\mathcal{A}_{r}$, $r<t$ which are incompatible with a $G\in 
\mathcal{B}_{t}$ are meaningless as noncausal for the observation at the
time $t$, they should be replaced by the causal histories $\mathcal{B}_{r}$, 
$r<t$ of the actual observables on $\mathfrak{h}$ which must be consistent
in the sense of compatibility of all $\mathcal{B}_{t}$. Thus the dynamical
formulation of the nondemolition principle of quantum causality and the
consistency of histories reads as 
\begin{equation*}
\mathcal{A}_{s}\subseteq \mathcal{B}_{r}^{\prime },\quad \mathcal{B}%
_{s}\subseteq \mathcal{B}_{r}^{\prime },\quad \forall r\leq s. 
\end{equation*}

These are the only possible conditions when the posterior states always
exist as results of inference (filtering and prediction) of future quantum
states upon the measurement results of the classical (i.e. commutative) past
of process of observation. The act of measurement transforms quantum
propensities into classical realities. As Lawrence Bragg, another Nobel
prize winner, once said, everything in the future is a wave, everything in
the past is a particle.

\subsubsection{The true Heisenberg principle}

The time continuous solution of the quantum measurement problem was
motivated by analogy with the classical stochastic filtering problem which
obtains the prediction of future for an unobservable dynamical process $%
x\left( t\right) $ by time-continuous measuring of another, observable
process $y\left( t\right) $. Such a problem was first considered by Wiener
and Kolmogorov who found its solution in the form of \ causal spectral
filter but only for the stationary Gaussian case. The differential solution
in the form of a stochastic filtering equation was then obtained by
Stratonovich \cite{Str66} in 1958 for an arbitrary Markovian pair $\left(
x,y\right) $. This was really a break through in the statistics of
stochastic processes which soon found many applications, in particular for
solving the problems of stochastic control under incomplete information (it
is possible that this was one of the reasons why the Russians were so
successful in launching the rockets to the Moon and other planets of the
Solar system in 60s).

If $X\left( t\right) $ is the unobservable process, a Heisenberg coordinate
process of a quantum particle, say, and $Y\left( t\right) $ is an observable
quantum process, describing the trajectories $y\left( t\right) $ of the
particle in a cloud chamber, say, why don't we find a filtering equation for
the a posterior expectation $q\left( t\right) $ of $X\left( t\right) $ in
the same way as we do it in the classical case if we know a history, i.e. a
particular trajectory $y\left( r\right) $ up to the time $t$? This problem
was first considered and solved for the case of quantum Markovian Gaussian
pair $\left( X,Y\right) $ corresponding to a quantum open linear system with
linear output channel, in particular for a quantum oscillator matched to a
quantum transmission line \cite{Be80}. By studying this example, the
nondemolition condition 
\begin{equation*}
\left[ X\left( s\right) ,Y\left( r\right) \right] =0,\quad \text{ }\left[
Y\left( s\right) ,Y\left( r\right) \right] =0\quad \forall r\leq s 
\end{equation*}
was first found, and this allowed to get the solution in the form of the
causal equation for $q\left( t\right) =\left\langle X\left( t\right)
\right\rangle _{y}$.

Let us describe this exact dynamical model of the causal nondemolition
measurement first in terms of quantum white noise analysis for a
one-dimensional quantum nonrelativistic particle of mass $m$ which is
conservative if not observed, in a potential field $\phi $. But we shall
assume that the particle is under indirect observation by measuring of its
Heisenberg position operator $X\left( t\right) $ with an additive random
error $e\left( t\right) :$%
\begin{equation*}
Y\left( t\right) =X\left( t\right) +e\left( t\right) . 
\end{equation*}
We take the simplest statistical model for the error process $e\left(
t\right) $, the white noise model (the worst, completely chaotic error),
assuming that it is a classical (i.e. commutative) Gaussian white noise
given by the first momenta 
\begin{equation*}
\left\langle e\left( t\right) \right\rangle =0,\quad \left\langle e\left(
s\right) e\left( r\right) \right\rangle =\sigma ^{2}\delta \left( s-r\right)
. 
\end{equation*}
The measurement process $Y\left( t\right) $ should be commutative,
satisfying the causal nondemolition condition with respect to the
noncommutative process $X\left( t\right) $ (and any other Heisenberg
operator-process of the particle), what can be achieved by perturbing the
particle Newton-Erenfest equation: 
\begin{equation*}
m\frac{\mathrm{d}^{2}}{\mathrm{d}t^{2}}X\left( t\right) +\nabla \phi \left(
X\left( t\right) \right) =f\left( t\right) . 
\end{equation*}
Here $f\left( t\right) $ is a Langevin force perturbing the dynamics due to
the measurement, which is assumed to be another classical (commutative)
white noise. 
\begin{equation*}
\left\langle f\left( t\right) \right\rangle =0,\quad \left\langle f\left(
s\right) f\left( r\right) \right\rangle =\tau ^{2}\delta \left( s-r\right) . 
\end{equation*}
In classical measurement and filtering theory the white noises $e\left(
t\right) ,f\left( t\right) $ are usually considered independent, and the
intensities $\sigma ^{2}$ and $\tau ^{2}$ can be arbitrary, even zeros,
corresponding to the ideal case of the direct unperturbing observation of
the particle trajectory $X\left( t\right) $. However in quantum theory
corresponding to the standard commutation relations 
\begin{equation*}
X\left( 0\right) =\hat{x},\quad \frac{\mathrm{d}}{\mathrm{d}t}X\left(
0\right) =\frac{1}{m}\hat{p},\quad \left[ \hat{x},\hat{p}\right] =i\hbar 
\hat{1} 
\end{equation*}
the particle trajectories do not exist, and it was always understood that
the measurement error $e\left( t\right) $ and perturbation force $f\left(
t\right) $ should satisfy a sort of uncertainty relation. This ''true
Heisenberg principle'' had never been mathematically formulated and proved
before the discovery \cite{Be80} of quantum causality and nondemolition
condition in the above form of commutativity of $X\left( s\right) $ and $%
Y\left( r\right) $ for $r\leq s$. As we showed first in the linear case \cite%
{Be80}, and later even in the most general case \cite{Be92a}, these
conditions are fulfilled if and only if $e\left( t\right) $ and $f\left(
t\right) $ satisfy the canonical commutation relations 
\begin{equation*}
\left[ e\left( r\right) ,e\left( s\right) \right] =0,\;\left[ e\left(
r\right) ,f\left( s\right) \right] =\frac{\hbar }{i}\delta \left( r-s\right)
,\;\left[ f\left( r\right) ,f\left( s\right) \right] =0. 
\end{equation*}
This proves that the pair $\left( e,f\right) $ must satisfy the uncertainty
relation $\sigma \tau \geq \hbar /2$, i.e. 
\begin{equation*}
\Delta e_{t}\Delta f_{t}\geq \hbar t/2\text{,} 
\end{equation*}
terms of the standard deviations of the integrated processes 
\begin{equation*}
e_{t}=\int_{0}^{t}e\left( r\right) \mathrm{d}r,\quad
f_{t}=\int_{0}^{t}f\left( s\right) \mathrm{d}s. 
\end{equation*}
This inequality constitutes the precise formulation of the true Heisenberg
principle for the square roots $\sigma $ and $\tau $ of the intensities of
error $e$ and perturbation $f$: they are inversely proportional with the
same coefficient proportionality $\hbar /2$ as for the to the pair $\left( 
\hat{x},\hat{p}\right) $. The canonical pair $\left( e,f\right) $ called
quantum white noise cannot be considered classically despite of the
possibility of the classical realizations of each process $e$ and $f$
separately due to the self-commutativity of the families $e$ and $f$. Thus,
a generalized matrix mechanics for the treatment of quantum open systems
under continuous nondemolition observation and the true Heisenberg principle
was invented exactly 20 years ago in \cite{Be80}.

\subsection{Consistent Histories and Filtering}

Schr\"{o}dinger believed that all quantum problems including the
interpretation of measurement should be formulated in continuous time in the
form of differential equations. He thought that the measurement problem
would have been resolved if quantum mechanics had been made consistent with
relativity theory and the time had been treated appropriately. However
Einstein and Heisenberg did not believe this, each for to his own reasons.
While Einstein thought that the probabilistic interpretation of quantum
mechanics was wrong, Heisenberg simply stated:-

\begin{quote}
\textit{Quantum mechanics itself, whatever its interpretation, does not
account for the transition from `possible to the actual'}
\end{quote}

Perhaps the closest to the truth was Bohr when he said that it `must be
possible so to describe the extraphysical process of the subjective
perception as\ if it were in reality in the physical world', extending the
reality beyond the closed quantum mechanical form by including a subjective
observer into a semiclassical world. He regarded the measurement apparatus,
or meter, as a semiclassical object which interacts with the world in a
quantum mechanical way but has only commuting observables - pointers. Thus
Bohr accepted that \textit{not all world is quantum mechanical, there is a
classical part of the physical world, and we belong partly to this classical
world.}

\subsubsection{Symbolic Stochastic Calculus}

In order to formulate the differential nondemolition causality condition and
to derive a filtering equation for the posterior states in the
time-continuous case we need quantum stochastic calculus.

The classical differential calculus for the infinitesimal increments 
\begin{equation*}
\mathrm{d}x=x\left( t+\mathrm{d}t\right) -x\left( t\right) 
\end{equation*}%
became generally accepted only after Newton gave a simple algebraic rule $%
\left( \mathrm{d}t\right) ^{2}=0$ for the formal computations of the
differentials $\mathrm{d}x$ for smooth trajectories $t\mapsto x\left(
t\right) $. In the complex plane $\mathbb{C}$ of phase space it can be
represented by a one-dimensional algebra $\mathfrak{a}=\mathbb{C}d_{t}$ of
the elements $a=\alpha d_{t}$ with involution $a^{\star }=\bar{\alpha}d_{t}$%
. Here 
\begin{equation*}
\text{$d_{t}$}=\left[ 
\begin{array}{ll}
0 & 1 \\ 
0 & 0%
\end{array}%
\right] =\frac{1}{2}\left( \sigma _{1}+i\sigma _{2}\right) 
\end{equation*}%
for $\mathrm{d}t$ is the nilpotent matrix, which can be regarded as
Hermitian $d_{t}^{\star }=d_{t}$ with respect to the Minkowski metrics $%
\left( \mathbf{z}|\mathbf{z}\right) =2\func{Re}z_{-}\bar{z}_{+}$ in $\mathbb{%
C}^{2}$.

This formal rule was generalized to non-smooth paths early in the last
century in order to include the calculus of forward differentials $\mathrm{d}%
w\simeq \left( \mathrm{d}t\right) ^{1/2}$ for continuous diffusions $w_{t}$
which have no derivative at any $t$, and the forward differentials $\mathrm{d%
}n\in \left\{ 0,1\right\} $ for left continuous counting trajectories $n_{t}$
which have zero derivative for almost all $t$ (except the points of
discontinuity where $\mathrm{d}n=1$). The first is usually done by adding
the rules 
\begin{equation*}
\left( \mathrm{d}w\right) ^{2}=\mathrm{d}t,\quad \mathrm{d}w\mathrm{d}t=0=%
\mathrm{d}t\mathrm{d}w 
\end{equation*}
in formal computations of continuous trajectories having the first order
forward differentials $\mathrm{d}x=\alpha \mathrm{d}t+\beta \mathrm{d}w$
with the diffusive part given by the increments of standard Brownian paths $%
w $. The second can be done by adding the rules 
\begin{equation*}
\left( \mathrm{d}n\right) ^{2}=\mathrm{d}n,\quad \mathrm{d}n\mathrm{d}t=0=%
\mathrm{d}t\mathrm{d}n 
\end{equation*}
in formal computations of left continuous and smooth for almost all $t$
trajectories having the forward differentials $\mathrm{d}x=\alpha \mathrm{d}%
t+\gamma \mathrm{d}m$ with jumping part given by the increments of standard
compensated Poisson paths $m_{t}=n_{t}-t$. These rules were developed by It%
\^{o} \cite{Ito51} into the form of a stochastic calculus.

The linear span of $\mathrm{d}t$ and $\mathrm{d}w$ forms the Wiener-It\^{o}
algebra $\mathfrak{b}=\mathbb{C}d_{t}+\mathbb{C}d_{w}$, while the linear
span of $\mathrm{d}t$ and $\mathrm{d}n$ forms the Poisson-It\^{o} algebra $%
\mathfrak{c}=\mathbb{C}d_{t}+\mathbb{C}d_{m}$, with the second order
nilpotent $d_{w}=d_{w}^{\star }$ and the idempotent $d_{m}=d_{m}^{\star }$.
They are represented together with $d_{t}$ by the triangular Hermitian
matrices 
\begin{equation*}
\text{$d_{t}$}=\left[ 
\begin{array}{lll}
0 & 0 & 1 \\ 
0 & 0 & 0 \\ 
0 & 0 & 0%
\end{array}%
\right] ,\quad d_{w}=\left[ 
\begin{array}{lll}
0 & 1 & 0 \\ 
0 & 0 & 1 \\ 
0 & 0 & 0%
\end{array}%
\right] ,\emph{\quad }d_{m}\mathbf{=}\left[ 
\begin{array}{lll}
0 & 1 & 0 \\ 
0 & 1 & 1 \\ 
0 & 0 & 0%
\end{array}%
\right] ,
\end{equation*}%
on the Minkowski space $\mathbb{C}^{3}$ with respect to the inner Minkowski
product $\left( \mathbf{z}|\mathbf{z}\right) =z_{-}z^{-}+z_{\circ }z^{\circ
}+z_{+}z^{+}$, where $z^{\mu }=\bar{z}_{-\mu }$, $-\left( -,\circ ,+\right)
=\left( +,\circ ,-\right) $.

Although both algebras $\mathfrak{b}$ and $\mathfrak{c}$ are commutative,
the matrix algebra $\mathfrak{a}$ generated by $\mathfrak{b}$ and $\mathfrak{%
c}$ on $\mathbb{C}^{3}$ is not: 
\begin{equation*}
d_{w}d_{m}=\left[ 
\begin{array}{lll}
0 & 1 & 1 \\ 
0 & 0 & 0 \\ 
0 & 0 & 0%
\end{array}%
\right] \neq \left[ 
\begin{array}{lll}
0 & 0 & 1 \\ 
0 & 0 & 1 \\ 
0 & 0 & 0%
\end{array}%
\right] =d_{m}d_{w}.
\end{equation*}%
The four-dimensional $\star $-algebra $\mathfrak{a}=\mathbb{C}d_{t}+\mathbb{C%
}e_{-}+\mathbb{C}e^{+}+\mathbb{C}e$ of triangular matrices with the
canonical basis 
\begin{equation*}
e_{-}=\left[ 
\begin{array}{lll}
0 & 1 & 0 \\ 
0 & 0 & 0 \\ 
0 & 0 & 0%
\end{array}%
\right] ,\,e^{+}\mathbf{=}\left[ 
\begin{array}{lll}
0 & 0 & 0 \\ 
0 & 0 & 1 \\ 
0 & 0 & 0%
\end{array}%
\right] ,\,e=\left[ 
\begin{array}{lll}
0 & 0 & 0 \\ 
0 & 1 & 0 \\ 
0 & 0 & 0%
\end{array}%
\right] ,
\end{equation*}%
given by the algebraic combinations 
\begin{equation*}
e_{-}=d_{w}d_{m}-\text{$d_{t}$},\;e^{+}=d_{m}d_{w}-\text{$d_{t}$}%
,\;e=e-d_{w}-\text{$d_{t}$}
\end{equation*}%
is the canonical representation of the differential $\star $-algebra for
one-dimensional vacuum noise in quantum stochastic calculus \cite{Be92, Be98}%
. It realizes the HP (Hudson -- Parthasarathy) table \cite{HuPa84} 
\begin{equation*}
\mathrm{d}\Lambda _{-}\mathrm{d}\Lambda ^{+}=\mathrm{d}t,\quad \mathrm{d}%
\Lambda _{-}\mathrm{d}\Lambda =\mathrm{d}\Lambda _{-},\quad \mathrm{d}%
\Lambda \mathrm{d}\Lambda ^{+}=\mathrm{d}\Lambda ^{+},\quad \left( \mathrm{d}%
\Lambda \right) ^{2}=\mathrm{d}\Lambda ,
\end{equation*}%
with zero products for all other pairs, for the multiplication of the
canonical counting $\mathrm{d}\Lambda =\Lambda \left( e\right) $, creation $%
\mathrm{d}\Lambda ^{+}=\Lambda \left( e^{+}\right) $, annihilation $\mathrm{d%
}\Lambda _{-}=\Lambda \left( e^{+}\right) $, and preservation $\mathrm{d}%
t=\Lambda \left( \text{$d_{t}$}\right) $ quantum stochastic integrators in
Fock space over $L^{2}\left( \mathbb{R}_{+}\right) $. As was proved in \cite%
{Be98}, any generalized It\^{o} algebra can be represented as a $\star $%
-subalgebra of a multi-dimensional quantum vacuum algebra 
\begin{equation*}
\mathrm{d}\Lambda _{\mu }^{\iota }\mathrm{d}\Lambda _{\kappa }^{\nu }=\delta
_{\kappa }^{\iota }\mathrm{d}\Lambda _{\mu }^{\nu },\quad \kappa ,\mu \in
\left\{ -,1,\ldots ,d\right\} ;\;\iota ,\nu \in \left\{ 1,\ldots
,d,+\right\} ,
\end{equation*}%
where $\mathrm{d}\Lambda _{-}^{+}=\mathrm{d}t$ and $d$ is the dimensionality
of quantum noise (could be infinite), similar to the representation of every
semi-classical system with a given state as a subsystem of quantum system
with a pure state. In particular, any real-valued process $y_{t}$ with zero
mean value $\left\langle y_{t}\right\rangle =0$ and independent increments
generating a two-dimensional It\^{o} algebra has the differential $\mathrm{d}%
y$ in the form of a commutative combination of $\mathrm{d}\Lambda ,\mathrm{d}%
\Lambda _{-},\mathrm{d}\Lambda ^{+}$. We shall call it standard if it has
the stationary increments with the standard variance $\left\langle
y_{t}^{2}\right\rangle =t$. In this case 
\begin{equation*}
y_{t}=\left( \Lambda ^{+}+\Lambda _{-}+\varepsilon \Lambda \right)
_{t}=\varepsilon m_{t}+\left( 1-\varepsilon \right) w_{t},
\end{equation*}%
where $\varepsilon \geq 0$ is defined by the equation $\left( \mathrm{d}%
y\right) ^{2}-\mathrm{d}t=\varepsilon \mathrm{d}y$.

\subsubsection{Stochastic decoherence equations}

The generalized wave mechanics which enables us to treat the quantum
processes of continual in time observation, or in other words, quantum
mechanics with trajectories $\omega =\left( y_{t}\right) $, was discovered
only quite recently, in \cite{Be88, Be89a, Be89b}. The basic idea of the
theory is to replace the deterministic unitary Schr\"{o}dinger propagation $%
\psi \mapsto \psi \left( t\right) $ by a linear causal stochastic one $\psi
\mapsto \chi \left( t,\omega \right) $ which is not necessarily unitary for
each history $\omega $, but unitary in the mean square sense with respect to
a standard probability measure $\mu \left( \mathrm{d}\omega \right) $ for $%
\mathrm{d}\omega $. Due to this positive measure 
\begin{equation*}
\Pr \left( t,\mathrm{d}\omega \right) =\left\| \chi \left( t,\omega \right)
\right\| ^{2}\mu \left( \mathrm{d}\omega \right) 
\end{equation*}
is normalized (if $\left\| \psi \right\| =1$) for each $t$, and is
interpreted as the probability measure for the histories $\omega _{t}=\left(
y_{r}\right) _{r<t}$.of an output stochastic process $y_{t}$. In the same
way as the abstract Schr\"{o}dinger equation can be derived from only
unitarity of propagation, the abstract decoherence wave equation can be
derived from the mean square unitarity in the form of a linear stochastic
differential equation. The reason that Bohr and Schr\"{o}dinger didn't
derive such equation despite their firm belief that the measurement process
can be described `as\ if it were in reality in the physical world' is that
the appropriate (stochastic) differential calculus had not been yet
developed early in that century. As Newton had to invent the differential
calculus in order to formulate the equations of classical dynamics, we had
to develop the quantum stochastic calculus for nondemolition processes \cite%
{Be88, Be92} in order to derive the generalized wave equation for quantum
dynamics with continual observation.

For the notational simplicity we shall consider here the one dimensional
case, $d=1$, the multi-dimensional case can be found elsewhere (e.g. in \cite%
{Be88, Be92}). The abstract stochastic wave equation can be written in this
case as \cite{Be89a, Be90b} 
\begin{equation*}
\mathrm{d}\chi \left( t\right) +K\chi \left( t\right) \mathrm{d}t=L\chi
\left( t\right) \mathrm{d}y_{t},\quad \chi \left( 0\right) =\psi 
\end{equation*}
where $y_{t}\left( \omega \right) $ is the output process, which is assumed
to be a martingale (e.g. the independent increment process with zero
expectations) representing a measurement noise with respect to the input
probability measure $\mu \left( \mathrm{d}\omega \right) =\Pr \left( 0,%
\mathrm{d}\omega \right) $ (but not with respect to the output probability
measure $\Pr \left( \infty ,\mathrm{d}\omega \right) $). The stochastic
process $\chi \left( t,\omega \right) $ normalized in the mean square sense
for each $t$, is a wave function $\chi \left( t\right) $ in an extended
Hilbert space, describing the process of continual decoherence of the
initial pure state $\rho \left( 0\right) =P_{\psi }$ into the mixture 
\begin{equation*}
\rho \left( t\right) =\int P_{\psi _{\omega }\left( t\right) }\Pr \left( t,%
\mathrm{d}\omega \right) 
\end{equation*}
of the posterior states corresponding to $\psi _{\omega }\left( t\right)
=\chi \left( t,\omega \right) /\left\| \chi \left( t,\omega \right) \right\| 
$. Assuming that the conditional expectation $\left\langle \mathrm{d}y_{t}%
\mathrm{d}y_{t}\right\rangle _{t}$ in 
\begin{eqnarray*}
\left\langle \mathrm{d}\left( \chi ^{\dagger }\chi \right) \right\rangle
_{t} &=&\left\langle \mathrm{d}\chi ^{\dagger }\mathrm{d}\chi +\chi
^{\dagger }\mathrm{d}\chi +\mathrm{d}\chi ^{\dagger }\chi \right\rangle _{t}
\\
&=&\chi ^{\dagger }\left( L^{\ast }\left\langle \mathrm{d}y_{t}\mathrm{d}%
y_{t}\right\rangle _{t}L-\left( K+K^{\ast }\right) \mathrm{d}t\right) \chi
\end{eqnarray*}
is $\mathrm{d}t$ (e.g. $\left( \mathrm{d}y_{t}\right) ^{2}=\mathrm{d}%
t+\varepsilon \mathrm{d}y_{t}$), the mean square normalization in its
differential form $\left\langle \mathrm{d}\left( \chi ^{\dagger }\chi
\right) \right\rangle _{t}=0$ can be expressed as $K+K^{\ast }=L^{\ast }L$,
\ i.e. 
\begin{equation*}
K=\frac{1}{2}L^{\ast }L+\frac{i}{\hbar }H, 
\end{equation*}
where $H=H^{\ast }$ is the Schr\"{o}dinger Hamiltonian such that this is the
Schr\"{o}dinger equation if $L=0$. One can also derive the corresponding
Master equation 
\begin{equation*}
\frac{\mathrm{d}}{\mathrm{d}t}\rho \left( t\right) +K\rho \left( t\right)
+\rho \left( t\right) K^{\ast }=L\rho \left( t\right) L^{\ast } 
\end{equation*}
for mixing the initially pure state $\rho \left( 0\right) =\psi \psi
^{\dagger }$, as well as a stochastic nonlinear wave equation for the
dynamical prediction of the posterior state vector $\psi _{\omega }\left(
t\right) $ which is the normalization of $\chi \left( t,\omega \right) $ at
each $\omega $.

Actually, there are two basic standard forms \cite{Be89b, Be90a} of such
stochastic wave equations, corresponding to two basic types of stochastic
integrators: the Brownian standard type, $\varepsilon =0$, $y_{t}\simeq
w_{t} $, and the Poisson standard type $\varepsilon =1$, $y_{t}\simeq
n_{t}-t $ with respect to the measure $\mu $. (For the most general form see 
\cite{Be97}.) To get these forms we shall assume that $y_{t}$ is standard
with respect to the input measure $\mu $, given by the multiplication table 
\begin{equation*}
\left( \mathrm{d}y\right) ^{2}=\mathrm{d}t+\nu ^{-1/2}\mathrm{d}y,\quad 
\mathrm{d}y\mathrm{d}t=0=\mathrm{d}t\mathrm{d}y, 
\end{equation*}
where $\nu >0$ is the intensity of the Poisson process $n_{t}=\nu
^{1/2}y_{t}+\nu t$ with respect to the standard measure $\mu $, and 
\begin{equation*}
L=\nu ^{1/2}(C-I),\quad H=E+i\frac{\nu }{2}\left( C-C^{\dagger }\right) . 
\end{equation*}
This corresponds to the stochastic decoherence equation for the counting
observation initially derived in the form \cite{Be89a, BaBe} 
\begin{equation*}
\mathrm{d}\chi \left( t\right) +\left( \frac{\nu }{2}\left( C^{\dagger
}C-I\right) +\frac{i}{\hbar }E\right) \chi \left( t\right) \mathrm{d}%
t=\left( C-I\right) \chi \left( t\right) \mathrm{d}n_{t}, 
\end{equation*}
where $C$ and $E$ are called collapse and energy operators respectively. The
nonlinear filtering equation for $\psi _{\omega }\left( t\right) $ in this
case has the form \cite{Be90a} 
\begin{equation*}
\mathrm{d}\psi _{\omega }+\left( \frac{\nu }{2}\left( C^{\dagger }C-\left\|
C\psi _{\omega }\right\| ^{2}\right) +\frac{i}{\hbar }E\right) \psi _{\omega
}\mathrm{d}t=\left( C/\left\| C\psi _{\omega }\right\| -I\right) \psi
_{\omega }\mathrm{d}n_{t}^{\psi _{\omega }}, 
\end{equation*}
where $n_{t}^{\psi }$ is the output counting process which is described by
the probability measure $\Pr \left( \infty ,\mathrm{d}\omega \right) $ with
the increment $\mathrm{d}n_{t}^{\psi }$ independent of $n_{t}^{\psi }$ under
the condition $\psi _{\omega }\left( t\right) =\psi $, having the
conditional expectation $\nu \left\| C\psi \right\| ^{2}\mathrm{d}t$. This
can be written also in the quasi-linear form \cite{Be89b, Be90a} 
\begin{equation*}
\mathrm{d}\psi \left( t\right) +\tilde{K}\left( t\right) \psi _{\omega
}\left( t\right) \mathrm{d}t=\tilde{L}\left( t\right) \psi \left( t\right) 
\mathrm{d}\tilde{y}_{t}^{\psi \left( t\right) }, 
\end{equation*}
where $\tilde{y}_{t}^{\psi }$ is the innovating martingale with respect to
the output measure which is described by the differential 
\begin{equation*}
\mathrm{d}\tilde{y}_{t}^{\psi }=\nu ^{-1/2}\left\| C\psi \right\| ^{-1}%
\mathrm{d}n_{t}^{\psi }-\nu ^{1/2}\left\| C\psi \right\| \mathrm{d}t 
\end{equation*}
with the initial $\tilde{y}_{0}^{\psi }=0$, the operator $\tilde{K}\left(
t\right) $ similar to $K$ has the form 
\begin{equation*}
\tilde{K}\left( t\right) =\frac{1}{2}\tilde{L}\left( t\right) ^{\ast }\tilde{%
L}\left( t\right) +\frac{i}{\hbar }\tilde{H}\left( t\right) , 
\end{equation*}
and $\tilde{H}\left( t\right) ,\tilde{L}\left( t\right) $ depend on $t$ and $%
\omega $ through the dependence on $\psi =\psi _{\omega }\left( t\right) $: 
\begin{equation*}
\tilde{L}=\nu ^{1/2}\left( C-\left\| C\psi \right\| \right) ,\quad \tilde{H}%
=E+i\frac{\nu }{2}\left( C-C^{\dagger }\right) \left\| C\psi \right\| . 
\end{equation*}

The last form of the nonlinear filtering equation admits the central limit $%
\nu \rightarrow \infty $ corresponding to the standard Wiener case $%
\varepsilon =0$ when $y_{t}=w_{t}$ with respect to the input Wiener measure $%
\mu $. If $L$ and $H$ do not depend on $\nu $, i.e. $C$ and $E$ depend on $%
\nu $ as 
\begin{equation*}
C=I+\nu ^{-1/2}L,\;E=H+\frac{\nu ^{1/2}}{2i}\left( L-L^{\dagger }\right) , 
\end{equation*}
then $\tilde{y}_{t}^{\psi }\rightarrow \tilde{y}_{t}$ as $\varepsilon
^{2}=\nu ^{-1}\rightarrow 0$, where $\tilde{y}_{t}$ with $\mathrm{d}\tilde{y}%
_{t}=\mathrm{d}y_{t}-2\func{Re}\left\langle \psi |L\psi \right\rangle 
\mathrm{d}t$ with respect to the output measure is another standard Wiener
process $\tilde{w}_{t}$ which does not depend on $\psi $. If $\left\| \psi
_{\omega }\left( t\right) \right\| =1$ (which follows from the initial
condition $\left\| \psi \right\| =1$), the operator-functions $\tilde{L}%
\left( t\right) $, $\tilde{H}\left( t\right) $ defining the nonlinear
filtering equation have the limits 
\begin{equation*}
\tilde{L}=L-\func{Re}\left\langle \psi |L\psi \right\rangle ,\quad \tilde{H}%
=H+\frac{i}{2}\left( L-L^{\dagger }\right) \func{Re}\left\langle \psi |L\psi
\right\rangle . 
\end{equation*}

\subsubsection{Quantum trajectories and filtering}

Let us give the exact model based on quantum stochastic calculus for a
quantum particle of mass $m$ in a potential $\phi $ under indirect
observation of its position given by the equation 
\begin{equation*}
\mathrm{d}Y_{t}=\left( 2\lambda \right) ^{1/2}X\left( t\right) \mathrm{d}t+%
\mathrm{d}y_{t}, 
\end{equation*}
where $y_{t}=w_{t}$ is the standard Wiener process. This model coincides
with the signal plus noise model $Y\left( t\right) =X\left( t\right)
+y\left( t\right) $ in terms of the generalized derivatives 
\begin{equation*}
y\left( t\right) =\left( 2\lambda \right) ^{-1/2}\mathrm{d}y_{t}/\mathrm{d}%
t\equiv e\left( t\right) , 
\end{equation*}
where $e\left( t\right) $ is quantum white noise of intensity $\sigma
^{2}=\left( 2\lambda \right) ^{-1}$. Here $\lambda >0$ is the accuracy of
the measurement of $X\left( t\right) $ with respect to the standard Wiener
process represented as $w_{t}=\left( \Lambda _{-}+\Lambda ^{+}\right) _{t}$
on the vacuum in Fock space. It was proved \ in \cite{Be88, Be92} that $%
Y_{t} $ is a commutative nondemolition process with respect to the
coordinate and momentum Heisenberg processes if they are perturbed by a
Langevian force $f\left( t\right) $ of intensity $\tau ^{2}=\lambda \hbar
^{2}/2$, the generalized derivative of $\left( \lambda /2\right) ^{1/2}f_{t}$%
, where $f_{t}=i\hbar \left( \Lambda _{-}-\Lambda ^{+}\right) _{t}$: 
\begin{equation*}
\mathrm{d}P\left( t\right) +\nabla \phi \left( X\left( t\right) \right) 
\mathrm{d}t=\left( \frac{\lambda }{2}\right) ^{1/2}\mathrm{d}f_{t},\quad
P\left( t\right) =m\frac{\mathrm{d}}{\mathrm{d}t}X\left( t\right) . 
\end{equation*}
Note that the quantum error process $e_{t}=w_{t}$ does not commute with the
perturbing quantum process $f_{t}$ in Fock space due to the multiplication
table 
\begin{equation*}
\mathrm{d}f\mathrm{d}w=i\hbar \mathrm{d}t,\quad \mathrm{d}w\mathrm{d}%
f=-i\hbar \mathrm{d}t\text{.} 
\end{equation*}
This corresponds to the canonical commutation relations for the normalized
derivatives $e\left( t\right) $ and $f\left( t\right) $, so that the true
Heisenberg principle is fulfilled at the boundary $\sigma \tau =\hbar /2$.
Thus our quantum stochastic model of nondemolition observation is the
minimal perturbation model for the given accuracy $\lambda $ of the
continual indirect measurement of the position $X\left( t\right) $.

The stochastic decoherence equation for this model with $y_{t}=w_{t}$, 
\begin{equation*}
L=\left( \frac{\lambda }{2}\right) ^{1/2}\hat{x},\quad H=\frac{\hat{p}^{2}}{%
2m}+\phi \left( \hat{x}\right) 
\end{equation*}
was derived and solved in \cite{Be88, Di88, Be89b, BeSt92} for the case of
linear and quadratic potential $\phi $. Here $\hat{x}$ is multiplication by $%
x$, and $\hat{p}=\frac{\hbar }{i}\frac{\mathrm{d}}{\mathrm{d}x}$.

The stochastic posterior nonlinear equation in this case is 
\begin{equation*}
\mathrm{d}\psi \left( t\right) +\left( \frac{i}{\hbar }H+\frac{\lambda }{4}%
\tilde{x}\left( t\right) ^{2}\right) \psi \left( t\right) \mathrm{d}t=\left( 
\frac{\lambda }{2}\right) ^{1/2}\tilde{x}\left( t\right) \psi \left(
t\right) \mathrm{d}\tilde{y}_{t}, 
\end{equation*}
where $\tilde{x}\left( t\right) =\hat{x}-q\left( t\right) \hat{1},\quad
q\left( t,\omega \right) =\mathrm{Tr}\left( \hat{x}P_{\psi _{w}\left(
t\right) }\right) $ is the posterior expectation (statistical prediction) of 
$\hat{x}$, and $\mathrm{d}\tilde{y}_{t}=\mathrm{d}y_{t}-\left( 2\lambda
\right) ^{1/2}q\left( t\right) \mathrm{d}t$ defines the innovating process $%
\tilde{y}_{t}$ which is equivalent to the standard Brownian motion $\tilde{w}%
_{t}$ with respect to the output (not input!) probability measure $\Pr
\left( \infty ,\mathrm{d}\omega \right) $.

Let us give the explicit solution of this stochastic wave equation for the
free particle ($\phi =0$) and the stationary Gaussian initial wave packet.
One can show \cite{ChSt92, Kol95} that the nondemolition observation of such
particle is described by filtering of quantum noise which results in the
continual collapse of any wave packet to the Gaussian stationary one
centered at the posterior expectation $q\left( t\right) =\left\langle \psi
\left( t\right) |\hat{x}\psi \left( t\right) \right\rangle $ with finite
dispersion $\left\| \left( q\left( t\right) -\hat{x}\right) \psi \left(
t\right) \right\| ^{2}\rightarrow \left( \hbar /2\lambda m\right) ^{1/2}$.
The center can be found from the linear Newton equation 
\begin{equation*}
\frac{\mathrm{d}^{2}}{\mathrm{d}t^{2}}z\left( t\right) +2\kappa \frac{%
\mathrm{d}}{\mathrm{d}t}z\left( t\right) +2\kappa ^{2}z\left( t\right)
=-g\left( t\right) , 
\end{equation*}
for the deviation process $z\left( t\right) =q\left( t\right) -y\left(
t\right) $, with $z\left( 0\right) =q_{0}-y\left( 0\right) $, $z^{\prime
}\left( 0\right) =v_{0}-y^{\prime }\left( 0\right) $. Here $\kappa =\left(
\lambda \hbar /2m\right) ^{1/2}$ is the decay rate which is also the
frequency of effective oscillations, $q_{0}=\left\langle \hat{x}%
\right\rangle $, $v_{0}=\left\langle \hat{p}/m\right\rangle $ are the
initial expectations, $y\left( t\right) $ is a (generalized) trajectory of
the observable process $Y\left( t\right) $, and $g\left( t\right) =y^{\prime
\prime }\left( t\right) $ is the effective gravitation for the particle in
the moving framework of $y\left( t\right) $. The continuous collapse $%
z\left( t\right) \rightarrow 0$ of the posterior expectation $q\left(
t\right) $ towards a linear trajectory $y\left( t\right) $ is illustrated in
the Appendix.

\subsubsection{A quantum message from the future}

Although Schr\"{o}dinger didn't derive the stochastic filtering equation for
the continuously decohering wave function $\chi \left( t\right) $,
describing the state of the semiclassical system including the observable
nondemolition process $y_{t}$ in continuous time in the same way as we did
it for his cat just in one step, he did envisage a possibility of how to get
it `if one introduces two symmetric systems of waves, which are traveling in
opposite directions; one of them presumably has something to do with the
known (or supposed to be known) state of the system at a later point in
time' \cite{Schr31}. This desire coincides with the ``transactional'' \
attempt of interpretation of quantum mechanics suggested in \cite{Crm86} on
the basis that the relativistic wave equation yields in the nonrelativistic
limit two Schr\"{o}dinger type equations, one of which is the time reversed
version of the usual equation: `The state vector $\psi $ of the quantum
mechanical formalism is a real physical wave with spatial extend and it is
identical with the initial ``offer wave'' of the transaction. The particle
(photon, electron, etc.) and the collapsed state vector are identical with
the completed transaction.' \ There was no mathematical proof of this
statement in \cite{Crm86}, and it is obviously not true for the
deterministic state vector $\psi \left( t\right) $ satisfying the
conventional Schr\"{o}dinger equation, but we are going to show that this
interpretation is true for the stochastic wave $\chi \left( t\right) $
satisfying our decoherence equation.

First let us note that the stochastic equation for the offer wave $\chi
\left( t\right) $ and the standard input probability measure $\mu $ can be
represented in Fock space as 
\begin{equation*}
\mathrm{d}\chi \left( t\right) +K\chi \left( t\right) \mathrm{d}t=L\mathrm{d}%
y_{t}\chi \left( t\right) ,\quad \chi \left( 0\right) =\psi \otimes \delta
_{0}, 
\end{equation*}
where $y_{t}=\Lambda ^{+}+\Lambda _{-}+\varepsilon \Lambda $. It coincides
on the noise vacuum state $\delta _{0}$ with the quantum stochastic Schr\"{o}%
dinger equation 
\begin{equation*}
\mathrm{d}\varphi \left( t\right) +K\varphi \left( t\right) \mathrm{d}%
t=\left( L\mathrm{d}\Lambda ^{+}-L^{\dagger }\mathrm{d}\Lambda _{-}\right)
\varphi \left( t\right) ,\quad \varphi \left( 0\right) =\psi \otimes \delta
_{0} 
\end{equation*}
corresponding to the generalized Heisenberg equation with the Langevin
force, $if_{t}=\hbar \left( \Lambda ^{+}-\Lambda _{-}\right) _{t}$, if $%
L^{\dagger }=L$. Indeed, due to adaptedness both $L\mathrm{d}y$ and $L%
\mathrm{d}\Lambda ^{+}-L^{\dagger }\mathrm{d}\Lambda _{-}$ act on the tensor
product states with future vacuum $\delta _{0}$ on which they have the same
action since $\Lambda _{-}\delta _{0}=0$, $\Lambda \delta _{0}=0$ (the
annihilation process $\Lambda _{-}$ is zero on the vacuum $\delta _{0}$, as
well as the number process $\Lambda $). Thus when extended from $\delta _{0}$
to any initial Fock vector $\varphi _{0}$, quantum stochastic evolution is
the HP unitary propagation \cite{HuPa84} which is a unitary cocycle on Fock
space over $L^{2}\left( \mathbb{R}_{+}\right) $ with respect to the free
time-shift evolution $\varphi \left( t,s\right) =\varphi \left( 0,s+t\right) 
$ in the Fock space. This free plain wave evolution in the half of space $%
s>0 $ in the extra dimension is the input, or offer wave evolution for our
three dimensional (or more?) world located at the boundary of $\mathbb{R}%
_{+} $. The single offer waves do not interact in the Fock space until they
reach the boundary $s=0$ where they produce the quantum jumps described by
the stochastic differential equation.

As has been recently shown in \cite{Be00}, by doubling the Fock space it is
possible to extend the cocycle to a unitary group evolution which will also
include the free propagation of the output waves in the opposite direction.
The conservative boundary condition corresponding to the interaction with
our world at the boundary, is including the creation, annihilation and
exchange of the input-output waves. The corresponding ``Schr\"{o}dinger''
boundary value problem is the second quantization of the Dirac wave equation
on the half line, with a boundary condition in Fock space which is
responsible for the stochastic interaction of quantum noise with our world
in the course of the transaction of the input-output waves. The
nondemolition continual observations are represented in this picture by the
measuring at the boundary of the arrival times and positions of the
particles corresponding to the quantized waves in Fock space with respect to
an ``offer state'' the input vacuum, dressed into the output wave. The
continual reduction process for our world wave function then is simply
represented as the decohering input wave function in the extended space,
which is filtered from the corresponding mixture of pure states by the
process of innovation of the initial knowledge during the continual
measurement. The result of this filtering gives the best possible prediction
of future states which is allowed by the quantum causality. As was shown on
the example of a free quantum particle under observation, the filtering
appear as a dissipation, oscillation and gravitation as a result of
nondemolition observation.

\section{Conclusion: The Greatest Form of Beauty}

My friend Robin Hudson wrote in his \ Lecture Notes on Quantum Theory:

\begin{quote}
\textit{Quantum theory is a beautiful mathematical theory. If only it didn't
have to mean something, to be interpreted.}
\end{quote}

Obviously here he used beautiful in the sense of simple: Everything that is
simple is indeed beautiful. However Nature is beautiful but not simple: we
live at the edge of two worlds, one is quantum, the other one is classical,
everything in the future is quantized waves, everything in the past is
trajectories of recorded particles.

In my philosophy I am a follower of those about whom Aristotle wrote in his
Metaphysica: `\textit{they fancied that the principles of mathematics are
the principles of all things'}, i.e. the things of Nature, and I agree that
`...\textit{these are the greatest forms of beauty'.}

\medskip

\textbf{Acknowledgment:}

I would like to acknowledge the help of Robin Hudson and some of my students
attending the lecture course on Modern Quantum Theory who were the first who
read and commented on these notes containing the answers on some of their
questions. The best source on history and drama of quantum theory is in the
biographies of the great inventors, Schr\"{o}dinger, Bohr and Heisenberg 
\cite{Schr, Bohr, Heis}, and on the conceptual development of this theory
before the rise of quantum probability --\ in \cite{Jamm}. An excellent
essay ``The quantum age begins'', as well as short biographies with posters
and famous quotations of all mathematicians and physicists mentioned here
can be found on mathematics website at St Andrews University --
http://www-history.mcs.st-and.ac.uk/history/, the use of which is
acknowledged.

\medskip

\textsc{Appendix}\textbf{: The continuous collapse of a free quantum
particle under the stationary nondemolition observation.}

The posterior position expectation $q\left( t\right) $ in the absence of
effective gravitation, $y^{\prime \prime }\left( t\right) =0$, collapses to
the registered linear trajectory $y\left( t\right) =ut-q$ with the rate $%
\kappa =\left( \lambda \hbar /2m\right) ^{1/2}$, remaining not collapsed, $%
q_{0}\left( t\right) =v_{0}t$ in the framework where $q_{0}=0$, only in the
classical limit $\hbar /m\rightarrow 0$ or absence of observation $\lambda
=0 $: 
\begin{equation*}
q_{0}\left( t\right) =v_{0}t,\quad q\left( t\right) =ut+e^{-\kappa t}\left(
q\cos \kappa t+\left( q+\kappa ^{-1}\left( v_{0}-u\right) \right) \sin
\kappa t\right) -q. 
\end{equation*}

\FRAME{itbpF}{17.0634cm}{10.0364cm}{0.424cm}{}{}{Plot}{\special{language
"Scientific Word";type "MAPLEPLOT";width 17.0634cm;height 10.0364cm;depth
0.424cm;display "PICT";plot_snapshots TRUE;mustRecompute FALSE;lastEngine
"MuPAD";xmin "0";xmax "6";xviewmin "0";xviewmax "6";yviewmin "-1";yviewmax
"2.0";viewset"XY";plottype 4;plottickdisable TRUE;num-y-ticks
1;labeloverrides 3;x-label "t";y-label "q(t)";numpoints 100;plotstyle
"patch";axesstyle "normal";xis \TEXUX{v58142};yis \TEXUX{y};var1name
\TEXUX{$\tau $};var2name \TEXUX{$y$};function \TEXUX{$5\tau $};linecolor
"maroon";linestyle 3;pointstyle "point";linethickness 1;lineAttributes
"Dots";var1range "0,6";num-x-gridlines 100;curveColor
"[flat::RGB:0x00800000]";curveStyle "Line";function \TEXUX{$0.5\tau
-1$};linecolor "blue";linestyle 3;pointstyle "point";linethickness
1;lineAttributes "Dots";var1range "0,6";num-x-gridlines 100;curveColor
"[flat::RGB:0x00000080]";curveStyle "Line";function \TEXUX{$e^{-\tau }\left(
6\sin \tau +\cos \tau \right) +0.5\tau -1$};linecolor "green";linestyle
1;pointstyle "point";linethickness 2;lineAttributes "Solid";var1range
"0,6";num-x-gridlines 50;curveColor "[flat::RGB:0x00008000]";curveStyle
"Line";rangeset"X";valid_file "T";tempfilename
'IRN8U102.wmf';tempfile-properties "XPR";}}

\end{document}